\documentclass{article}
\usepackage[title,titletoc,toc]{appendix}
\usepackage{amsmath}
\usepackage{amsfonts}
\usepackage[latin1]{inputenc}
\usepackage{graphicx}
\usepackage{rotating}
\usepackage{float}
\newcommand{\Ref}[1]{(\ref{#1})}
\newcommand{\Frac}[2]{\frac{\textstyle #1}
                          {\textstyle #2}}

\newcommand{\BigI}{{\rm 1\kern-.17em 1}}

\newcommand{\BigR}{{\rm I\kern-.17em R}}
\newcommand{\BigN}{{\rm I\kern-.17em N}}
\newcommand{\BigC}{{\rm \kern.24em \vrule width.02em height1.4ex depth-.05ex \kern-.26em C}}

\newcommand{\di}{\displaystyle}

\newcommand{\R}{\mathbb  R}
\newcommand{\C}{\mathbb  C}
\newcommand{\N}{\mathbb  N}
\newcommand{\T}{\mathbb  T}
\newcommand{\Z}{\mathbb  Z}
\newcommand{\1}{\mathbb  I}
\newcommand{\la}{\langle}
\newcommand{\ra}{\rangle}

\newtheorem{lemma}{Lemma}[section]
\newtheorem{theorem}[lemma]{Theorem}
\newtheorem{proposition}[lemma]{Proposition}

\newtheorem{remark}[lemma]{Remark}

\newtheorem{definition}[lemma]{Definition}

\def\noi{\noindent}

\def\beq{\begin{equation}}   \def\eeq{\end{equation}}
\def\bea{\begin{eqnarray}}  \def\eea{\end{eqnarray}}

\def\noi{\noindent}

\newcommand\mysection{\setcounter{equation}{0}\section}
\renewcommand{\theequation}{\thesection.\arabic{equation}}
\newcounter{hran} \renewcommand{\thehran}{\thesection.\arabic{hran}}

\def\bmini{\setcounter{hran}{\value{equation}}
    \refstepcounter{hran}\setcounter{equation}{0}
    \renewcommand{\theequation}{\thehran\alph{equation}}\begin{eqnarray}}

\def\bminiG#1{\setcounter{hran}{\value{equation}}
\refstepcounter{hran}\setcounter{equation}{-1}
\renewcommand{\theequation}{\thehran\alph{equation}}
\refstepcounter{equation}\label{#1}\begin{eqnarray}}
 
\begin{document}

\title { Numerical approaches for some Nonlinear  Eigenvalue Problems}
\author{ Fatima Aboud \\
Mathematics Department, College of Science, University of Diyala, Iraq \\
{\it \small Fatima.Aboud@sciences.uodiyala.edu.iq} \\
Fran\c{c}ois Jauberteau, Guy Moebs \&  Didier Robert \\
Laboratoire de Math\'ematiques Jean Leray, CNRS-UMR 6629, \\
 Universit\'e de Nantes, France \\
{\it \small francois.jauberteau@univ-nantes.fr,\; 
guy.moebs@univ-nantes.fr,\; didier.robert@univ-nantes.fr}
}
\vskip 1 truecm
\date{}
\maketitle

\begin{abstract}
In this article we are interested for the numerical study of nonlinear eigenvalue problems. We begin with a review of theoretical results obtained by functional analysis methods, especially for the Schr\"odinger pencils. Some recall are given for the pseudospectra. Then we present the numerical methods and results obtained for eigenvalues computation with spectral methods and finite difference discretization, in infinite or bounded domains. Comparison with theoretical results is done. The main difficulty here is that we have to compute eigenvalues of strongly non-self-adjoint operators 
  which are very unstable.  
  \end{abstract}

\noindent {\bf Keywords~:}  nonlinear eigenvalue problems, spectra, pseudospectra, finite difference methods, Galerkin spectral method, Hermite functions.

\pagestyle{myheadings}

\mysection{Introduction }
We are interested  here  in  equations   like $L(\lambda)u=0$  where  $L(\lambda)$  is  a   linear  operator on some linear space ${\cal E}$, depending on
 a  complex parameter $\lambda$.  When $L(\lambda)  = L_0 -\lambda\1$, this is the usual  eigenvalue problem : find  $\lambda\in\C$  and $u\in{\cal E}$,
  $u\neq 0$  such that $L(\lambda)u=0$. \\
In many  applications, in particular for  dissipative problems in mechanics, it is necessary to consider more  general  dependance in  the complex parameter $\lambda$.
  A particular  interesting case  is  a  quadratic  dependence : $L(\lambda) =   \lambda^2L_2 + \lambda L_1 + L_0$.
   We shall say  that  $L(\lambda)$ is a quadratic pencil.\\
 Let us consider the second order differential equation
 \beq\label{Deq1}
 \frac{d^2u}{dt^2}L_2+ \frac{du}{dt}L_1+ uL_0 = 0
 \eeq
Equation (\ref{Deq1})  is a model in mechanics for small oscillations of a continuum
   system in the presence of an impedance force \cite{krla}.\\
  Now looking for stationary solutions of  (\ref{Deq1}),   $u(t) = u_0{\rm e}^{\lambda t}$,
   we have the following equation
   \beq\label{Deq2}
   (\lambda^2L_2 + \lambda L_1 + L_0)u_0= 0
   \eeq
   So equation (\ref{Deq2}) is a non linear eigenvalue problem in the spectral parameter $\lambda\in \mathbb C$.\\
   The operator  $L_1$ represent  a  damping  term  as we see  in the  following  simple example.\\
   Let us consider  the  perturbed  wave equation
   \beq\label{W1D}
   \frac{\partial^2}{\partial t^2}u -\frac{\partial^2}{\partial^2x}u -2a \frac{\partial}{\partial t}u = 0
   \eeq
   where $t\in\R$  and $x\in\T := \R/2\pi\Z$. The  damping  term  $a <0$  is  here  constant. 
   So  we  have   to   solve   (\ref{W1D}) with   periodical boundary conditions.\\
   The  stationary  problem  is reduced  to the equation
   $$
   \lambda^2  + k^2 -2a\lambda = 0,\;\; k\in\Z
   $$
   Then we have for $k^2\geq a^2$ the  damped solutions  of (\ref{W1D}) : 
   $$
   u_k(t,x) = \exp\left((a+i\sqrt{k^2-a^2})t+ ikx\right)
   $$
   When  $a$ is a function  of $x$  we have  no  explicit formula  so   we need  numerical approximations
    to compute  the  damping modes. It is  the main goal  of this work, in particular concerning the Schr\"odinger   pencil
     $L_{V,a}(\lambda)$.

   We say that $\lambda$ is a non linear eigenvalue if there  exists $u_0\neq 0$  satisfying (\ref{Deq2}).\\
   Such  generalized   eigenvalue problems  have appeared   in a completely different  way. 
   The question  was to decide  if  a class  of  P.D.E  with analytic  coefficients  preserves  or not  the  analyticity property. To be more  explicit, let us consider  a P.D.E :  $Pu =f$. Assume that  $f$  is analytic 
   in some  open set $\Omega$, is-it  true  that $u$  is analytic in  $\Omega$ ? This is true  for  elliptic operators.\\
   For     some example, this question can be reduce to the    following  (see \cite{he2}  for more details): \\
   Does there  exist $\lambda\in\C, 0\neq u\in{\cal S}(\R)$  such  that 
    \beq\label{Deq3}
        \left(-\frac{d^2}{dx^2} + (x^2 -\lambda)^2\right)u = 0 \ \ ?
        \eeq

    Existence of  non null solutions for (\ref{Deq2})  and  (\ref{Deq3})  is a non trivial problem. 
  For  (\ref{Deq3}) it was  solved in \cite{phro} where it is proved that the generalized eigenfunctions span  the Hilbert  space $L^2(\R)$.\\
On the other side   we can prove that the  equation 
$$
    \left(-\frac{d^2}{dx^2} + (x -\lambda)^2\right)u = 0
$$
  has only the  trivial   solution  $u\equiv 0$ in $L^2(\R)$, $\forall \lambda\in\C$.
  
  Our aim  in this work is to   present   several   numerical approaches concerning  this  kind of non linear eigenvalue problems. \\
  For simplicity we  only  consider  quadratic  pencil such that $L_0 =\1$. 
  We can  reduce to this case  if $L_0$ or $L_2$ are   invertible  in the  linear space ${\cal E}$.\\
  To  every  quadratic pencil $L(\lambda)$  we can associate  a linear operator $ {\cal A}_L$ in ${\cal E}\times{\cal E}$ 
  such  that $\lambda$  is  a  non linear   eigenvalue  for $L$ if and only if $\lambda$  is a usual eigenvalue for  ${\cal A}_L$.  \\
  $ {\cal A}_L$  is called  a linearization  of  $L(\lambda)$. 
  It is  easy  to see that  we can choose
$$
                      {\mathcal  A}_L = \left(
                            \begin{array}{cc}
                            0 & \1\\
                            -L_0&-L_1\\
                            \end{array}
                            \right)
$$
                   So   non-linear eigenvalue  problems (for polynomial  operator  pencils)   can be reduced  to usual  eigenvalue problems
                    but   it is useful   to take care  of their  particular  structure. There  exist   infinitely many  linearizations.
                    
                    We are  mainly interested  here in  the  multidimensional  case   called   Schr\"odinger  pencils:
                    $$
                    L_{V,a}(\lambda) = -\triangle + V -2a\lambda  +\lambda^2
                    $$
                    in the  Hilbert  space $L^2(\R^d)$. $V$ and $a$  are  smooth  real functions  on $\R^d$  such  that 
                    $\di{\lim_{\vert x\vert \rightarrow +\infty}V(x) = +\infty}$  and $\vert a\vert\leq\sqrt V$.\\
                    The  main questions  we want to discuss  is the location  in  the  complex plane $\C$  of the
                    eigenvalues of $L_{V,a}$.   In the  first  part of this work  we shall recall  some  known theoretical
                      results  and in the  second part  we shall discuss several  numerical approaches for  the  computation
                       of the  eigenvalues  of $L_{V,a}$.  We shall see that  accurate  theoretical  results  on the location
                       of the  eigenvalues have been obtained  for 1D pencils $L_{V,a}$  but in the multidimensional
                        case very  few results  are known on the eigenvalues  of $L_{V,a}$   when $a$  is of the same
                     order  of $\sqrt V$.\\
                    
In Section 2, we present a review of theoretical results obtained by functional analysis methods. In Section 3, we give more results  for Schr\"odinger pencils. In Section 4 we recall some results on pseudospectra. In Section 5 we present the numerical methods and results obtained for eigenvalues computation with spectral methods and finite difference discretization, in infinite and bounded domain. Comparison with theoretical results is done. Then in Section 6 we give conclusions and open problems. 
                    
  \section{  A review  of  theoretical results obtained  by functional analysis methods}
               Most  of these results  was obtained by the Russian school between 
               1917 and 1970. For  more details we refer to the book \cite{ma}.
               
                Let us consider the quadratic family of operators $L(\lambda)=L_0+\lambda L_1 +\lambda^2$
       where $L_0$, $L_1$ are operators in an Hilbert space ${\mathcal H}$.\\
       If ${\mathcal H}$  if  of  dimension $N<+\infty$   the eigenvalues  are  the  solutions of
        the   polynomial equation $\det(L(\lambda)) = 0$. When $N$ is large this could be a difficult  problem at least  for numerical  computations.\\
         In applications  involving PDE,   ${\mathcal H}$  is a 
        $L^2$ space or a  Sobolev space, which is infinite  dimensional  and there is no explicit  equation  for the  generalized eigenvalues. Moreover, as we shall see later,  the  non  linear  eigenvalue problem  is equivalent to a linear
        eigenvalue problem which, in general, is  non self-adjoint hence unstable.

       $L_0$ is assumed to be self-adjoint, positive, with a domain $D(L_0)$
       and $L_1$ is $\sqrt{L_0}$-bounded. Moreover $L_0^{-1/2}$ is in a Schatten class
       ${\mathcal C}^p({\mathcal H})$\footnote{Recall that  a compact operator $A$    in an Hilbert  space  is in the Schatten
        class ${\mathcal C}^p$ if the  series $s_j(A)$  of  the eigenvalues   of $\sqrt{A^*A}$ satisfies $\sum s_j(A)^p < +\infty$},    for some real $p>0$.   \\
          The following results are well known.
          \begin{theorem}\label{gene}
          $L(\lambda)$ is a family of closed operators in ${\mathcal H}$. \\
          $\lambda\mapsto L^{-1}(\lambda)$ is meromorphic in the complex plane.\\
          The poles  $\lambda_j$ of $L^{-1}(\lambda)$,   with multiplicity $m_j$, 
          co\"{\i}ncide with the eigenvalues with the same multiplicities, of the matrix
          operator  ${\mathcal A}_L$ in the Hilbert space ${\mathcal H}\times D(L_0^{1/2})$,
           with domain $D({\mathcal A}_L)=D(L_0)\times D(L_0^{1/2})$
           where 
$$
                           {\mathcal  A}_L = \left(
                            \begin{array}{cc}
                            0 & \1\\
                            -L_0&-L_1\\
                            \end{array}
                            \right)
$$
                            \end{theorem}   
                          Assuming  that $V(x) \geq C\vert x\vert^{2m}$ and
                           $\vert a(x)\vert \leq C\sqrt{V(x)}$, $C>0$,  then   the Schr\"odinger  pencil $L_{V,a}(\lambda)$  satisfies the  above theorem  for $p>\frac{d(m+1)}{2m}$.\\
                   
                   If $L_0$  is positive and non degenerate       we have the  symmetric linearization
$$
                    {\mathcal  A}_{SL} = \left(
                            \begin{array}{cc}
                            0 & \sqrt{L_0}\\
                            -\sqrt{L_0}&-L_1\\
                            \end{array}
                            \right)
$$

                              Let us denote Sp[$L$] the eigenvalues of   ${\mathcal A}_L$ (which co\"{\i}ncide with the poles of $L^{-1}(z)$).
             \begin{remark}
             It may happens that ${\rm Sp}[L]$is empty  (example : $L(\lambda) = -\frac{d^2}{dx^2} + (x -\lambda)^2$).
             \end{remark}
             Let us remark that if $L_1=0$   then $\lambda\in {\rm Sp}[L]$ if and only if $-\lambda^2$  is in the spectrum of $L_0$. So if $L_0$ has  a point spectrum  then
            ${\rm Sp}[L]$ is a subset  of the imaginary  axis.

            We shall see now that when $L_1$  is  strictly smaller  than   $\sqrt L_0$  then the   eigenvalues  are asymptotically close  to 
            the  imaginary axis  and the  generalized  eigenvectors is a dense set in the  Hilbert space. 
             When  $L_1$  has the  same power  of      $\sqrt L_0$  
             it may happens that there  is  no eigenvector at all for $L(\lambda)$.
             
             If $\lambda_0\in{\rm Sp}[L]$  we denote by ${\cal E}_L(\lambda_0)$  the linear  space of the solutions 
             $\{u_0,u_1,\cdots, u_k,\cdots \}$ of the equations
$$
             L(\lambda_0)u_0 =0,\;\; L(\lambda)u_1+ L^\prime(\lambda_0)u = 0 ,\\
             L(\lambda_0)u_{k+2} + L^\prime(\lambda_0)u_{k+1} + \frac{1}{2}L^{\prime\prime}(\lambda_0)u_k=0,\; k\geq 0
$$
     The  dimension of  ${\cal E}_L(\lambda_0)$  is  the multiplicity of $\lambda_0$ (for details see \cite{phro}).\\
     Assume  that $L_0, L_1$  are  self-adjoint,  $L_0$  is positive non degenerate  and that 
     there  exist $\kappa \geq 0$  and $\delta \geq 0$ such that $L_1L_0^{\delta-1/2}$  is a bounded  operator on 
     ${\cal H}$  and
     $$
     \Vert L_1L_0^{\delta-1/2}\Vert \leq \kappa 
     $$
     Assume that $L_0^{-1}$ is in the  Schatten class $C^p$, $p\geq 1$.
     \begin{theorem}\label{thm2}
     If $0<\delta \leq 1/2$   then the  spectra of $L$  is the  domain
     $$
     \Omega_\delta =  D_R\cup \{\lambda\in\C,\; \vert\Re\lambda\vert \leq \kappa\vert\lambda\vert^{1-2\delta}\}
     $$ 
      and $\oplus_{\lambda\in{\rm Sp}[L]}{\cal E}_L(\lambda)$  is  dense  in ${\cal H}$.\\
     If $\delta =0$ and if 
     $$
     \vert\frac{\pi}{2} - \arccos\kappa\vert \leq \frac{\pi}{2p}
     $$
     then  $\oplus_{\lambda\in{\rm Sp}[L]}{\cal E}_L(\lambda)$  is also  dense  in ${\cal H}$. 
     \end{theorem}
  For $\delta>0$ we get that the  eigenvalues  are localized in a vertical parabolic domain in the imaginary direction. For $\delta =0$   end $\kappa$  small  the  eigenvalues are localized  in a  small sector  around  the imaginary axis.
  Notice that for $\kappa$ of order 1 the above  theorem does not  give any information on the
   location of ${\rm Sp}[L]$; we only know that it  is a  discrete and  infinite   subset of $\C$.\\
   
\noindent {\bf A sketch of proof of Theorem (\ref{thm2})}\\
      The idea  is to  consider $L(\lambda)$  as a perturbation  of $L_0+\lambda^2$.
      We know that $L_0+\lambda^2$ has a spectrum  in $i\R$  because  $L_0$  is self-adjoint.\\
      We  have
$$
      L(\lambda) = \left(\1+\lambda L_1)(L_0+\lambda^2)^{-1}\right)(L_0+\lambda^2)
$$
      So if  $\lambda\notin i\R$   then $L(\lambda)$  is invertible  if and only  $(\1+\lambda L_1)(L_0+\lambda^2)^{-1}$ is invertible. To check  this  property it is enough  to choose  $\lambda$ such that 
      $\Vert\lambda L_1(L_0+\lambda^2)^{-1}\Vert <1$.  
      $\square$
      
      \noi
  Moreover  If $L_1$  has  a  sign   we have easily
   \begin{proposition}
   If $L_1\geq 0$   then ${\rm Sp}[L]\subseteq\{\lambda\in\C, \Re\lambda\leq 0\}$.\\
   If If $L_1\leq0$   then ${\rm Sp}[L]\subseteq\{\lambda\in\C, \Re\lambda\geq 0\}$.
   \end{proposition}
{\bf Proof}.  If $L(\lambda)u = 0 $  then $\langle u, L(\lambda)u\rangle = 0$.
 Taking the imaginary  part of this  equality we get the  proposition. $\square$\\
  The above result  applies for example  to 
  $$
  L(\lambda)  =  -\frac{d^2}{dx^2} + x^6 +\alpha x^2\lambda + \lambda^2
  $$
  For  this example  we have $\delta = \frac{5}{6}$  hence  the spectra is localized  inside  the  parabolic region \\
   $\{\lambda\in\C,\; \vert\Im\lambda\vert \geq C\vert\Re\lambda\vert^{5/2}\}$.

       For Schr\"odinger pencils $L_{V,a}$  we can say  more.
 \section{More results  for Schr\"odinger pencils}
 Let us recall our  definition  of Schr\"odinger  pencils:
 $  L_{V,a}(\lambda) = -\triangle + V -2a\lambda  +\lambda^2$.\\

 In all this article we assume that  the  pair of functions$(V,a)$ satisfies  the  following  technical  conditions.
 We do not try here  to discuss the optimality of this conditions.\\
  $[{\rm cond}(V,a)]$.  $V, a$  are  smooth  $C^\infty$  functions  on $\R^d$. There  exists $k>0$ such that
\bea
  \vert\partial_x^\alpha V(x)\vert \leq C_\alpha\la x\ra^{k-\vert\alpha\vert},\;\;
   \vert\partial_x^\alpha a(x)\vert \leq C_\alpha\la x\ra^{k/2-\vert\alpha\vert},\\
  \vert a(x)\vert \leq\sqrt{V(x)},\;\; V(x) \geq 0,\; V(x)\geq c\la x\ra^k,\; {\rm for}\; \vert x\vert\geq 1
    \eea
with $C_\alpha >0$ and $c >0$. 
    Under these  conditions we know that $L_0 =-\triangle +V$  is an unbounded  self-adjoint  operator
     in $L^2(\R^d)$  and for  every $\lambda\in\C$  $L_{V,a}(\lambda)$ is  a closed and  Fredholm  operator  with domain the
     following  weighted  Sobolev space: 
     ${\cal H}_V=\{u\in L^2(\R^d),\; \triangle u\in  L^2(\R^d), Vu\in  L^2(\R^d)\}$.  Moreover the  set  
     ${\rm Sp}[L]$ of eigenvalues  of $L_{V,a}$  is a discrete set (empty or not), each  eigenvalue having
       a finite  multiplicity and the only  possible accumulation point in the complex plane is $\infty$. \\
       Notice  that  $\lambda$  is an eigenvalue  then its complex conjugate $\bar\lambda$ is also  an eigenvalue.
\begin{proposition}
Assume  that  $(V, a)$  satisfies $[{\rm cond}(V,a)]$ and that  $a\ \leq 0$,  $a(x^0)<0$ 
 for some  $x^0\in\R^d$. Then  ${\rm Sp}[L]  $  is in the open sector
 $\{\lambda\in\C; \,\; \Im(\lambda)>0,\; \Re(\lambda)\neq 0\}$.
 \end{proposition}
{\bf  Proof} Let $u\in L^2(\R^d)$, $u\neq 0$ such that $L_{V,a}(\lambda)u=0$. Set $\lambda = r+is$.
 We know that $r\geq 0$.  Assume that $r=0$. Reasoning by contradiction we first prove that $s=0$.
 If $s\neq 0$  that  we  get  that $\int_{\R^d}a(x)\vert u\vert^2(x)dx = $  hence $u$  vanishes  in an  
 non empty ope set of $\R^d$  and applying the  uniqueness  Calderon theorem for second order
  elliptic equation  we get $u=0$ on $\R^d$  and a contradiction.\\
  If $s=0$        we get 
  $$
  (-\triangle + V-2ra + r^2)u = 0
  $$
  and  $\int_{\R^d} (V(x)-2ra(x) + r^2))\vert u(x)\vert^2dx = 0$. Using that  $V\leq a^2$
   we have $ \int_{\R^d} (r-a)^2\vert u(x)\vert^2dx = 0$. So again we get that $u$ vanishes
    on a  non empty open  set and a contradiction like above. $\square$
    
    Let us remark that  the general results given   in Theorem \ref{thm2}  apply if there exists
     $\delta\geq 0$ such that $\vert a\vert(x)\leq C V(x)^{1/2-\delta} $ or
     $\vert a(x)\vert\leq \kappa V(x)^{1/2} $ with $\kappa$ small  enough.
     
     For 1D  Sch\"odinger  pencils accurate  results  were  obtained by M. Christ \cite{chr1,chr2}  et by \cite{ccy}.  Let us recall  here   some  of their  results.
     They  consider  the pencils
      $$
       L_k(\lambda) = -\frac{d^2}{dx^2} + (x^k -\lambda)^2
       $$  
       with $k\in\N$. Here  we shall only  consider $k$ even. The above assumptions are  satisfied.\\
       \begin{proposition}[M. Christ \cite{chr3}]
       For  every $k\geq 2$, $k$ even,  the set  ${\rm Sp}[L_k]$  is included  in  the two  sectors
        $\{\lambda\in\C,\; \vert\arg(\lambda)\vert \geq \frac{k\pi}{2(k+1)}\}$.
     \end{proposition}
     The second result  say  that  the  eigenvalues  of  large  modulus  are  close  to
      the  lines $\{\lambda\in\C,\; \vert\arg(\lambda)\vert \geq \frac{k\pi}{2(k+1)}$.
      \begin{theorem}[Y. Ching-Chau, \cite{ccy}, Theorem 1]
      Let $\{\lambda_n\}_{n\in\N}$  be the  set ${\rm Sp}[L_k]$  such that
       $\vert\lambda_1\vert < \vert\lambda_2\vert<\cdots<\vert\lambda_n\vert<\vert\lambda_{n+1}\vert<\cdots$.\\    Then  we   have
        for $n\rightarrow +\infty$, 
       \beq\label{asympt}
       \lambda_n= \left(\frac{\pm(n+\frac{1}{2})\pi i - \log(2)}{\frac{2k}{k+1}}\right)^{\frac{k}{k+1}} + 
       O\left(\frac{1}{\sqrt n}\right)
       \eeq
        \end{theorem}
        This result was  proved  using   ODE  methods  in the  complex plane.\\
      By  an elementary  computation of  the argument  for  the complex number in the r.h.s of (\ref{asympt})    we can see that
      $\vert\arg(\lambda_n)\vert$ is  close  to $\frac{k\pi}{2(k+1)}$  when $n\rightarrow +\infty$.
     We also  have the following result
     \begin{theorem}[\cite{phro, ab}]
    The  linear  space  span  by  the  generalized   eigenfunctions associated  with the eigenvalues
     $\{\lambda_n\}$   is dense in $L^2(\R)$.
     \end{theorem}
     In $\cite{phro}$  the proof   was  given for $L_2(\lambda)$  and for  $L_k(\lambda$, $k>2$, even in \cite{ab}, \cite{abro}.

     In the  following result  we shall  see that the spectral  set  ${\rm Sp}(L_k)$  is very  unstable  under
      perturbations. M. Christ \cite{chr2}  has consider  the  following  model:
$$
      L_P^\#(\lambda) = (P-\lambda +\frac{d}{dx})(P-\lambda -\frac{d}{dx})
$$
        We also  have  $ L_P^\#(\lambda)  = -\frac{d^2}{dx^2} + (P-\lambda)^2 + P^\prime$;
         where  $P$ is a polynomial. Assume   that  the  degree $k$ of $P$  is even,
          $P(x) = x^{k} + a_{k-1}x^{k-1}+\cdots+a_1x+a_0$.
          \begin{proposition}\label{ce}
          We have  ${\rm Sp}[L_P^\#]=\emptyset$. In other words for every $\lambda\in\C$,
           the equation  $ L_P^\#(\lambda) u=0$ has only the  trivial solution $ u\equiv 0$ \footnote{it is known that every
            solution  in  $L^2(\R)$ of   $ L_P^\#(\lambda) u=0$ is in the  Schwartz  space
           ${\cal S}(\R)$ (see \cite{phro})}  in the  Schwartz  space
           ${\cal S}(\R)$. 
      \end{proposition}                                                                                                                                                                                                                                                                                                                                                                                                                                                                                                                                                                                                                                                                                                                                                                                                                                                                                                                                                                                                                                                                                                                                                                                                                                                                                                                                                                                                                                                                                                                                                                                                                                                                                                                                                                                                                                                                                                                       
{\bf A sketch of proof of Proposition (\ref{ce})}\\
 We have 
 $$
  L_P^\#(\lambda) = (P-\lambda +\frac{d}{dx})(P-\lambda -\frac{d}{dx})
  $$
So, we have to solve the two equations
\bea
(P-\lambda +\frac{d}{dx})v &=0 \\
(P-\lambda -\frac{d}{dx})u &=v
\eea
Set $Q(x) = \int_0^x(P(s)-\lambda)ds$  and using standard ODE methods  we get that
\beq\label{sol}
u(x) = C_1 {\rm e}^{Q(x)} + C_2 {\rm e}^{Q(x)}\int_x^{+\infty}{\rm e}^{_2Q(s)}ds 
\eeq
 where $C_1, C_2$  are constants.  If $u$ is in the Schwartz space then $u$  is in particular  bounded
  but (\ref{sol})  shows that this is possible  only if $C_1=C_2=0$.$\square$


\section{ Pseudospectra for linear pencils}

As we have  seen above    the    eigenvalues  of Schr\"odinger  pencils  are  very  unstable.  As  propose some  times ago by Thefthen \cite{thr1}   it  is  useful  to replace the  spectra of non-self adjoint operators   by  something   more stable which is called  the  pseudospectra.

\subsection{A short review}

Let $A$  be  closed operator in the Hilbert  space ${\cal H}$  with domain $D(A)$  dense in ${\cal  H}$. 
Recall that $D(A)$ is an Hilbert space for  the graph norm 
$\Vert u\Vert_{D(A)} = \sqrt{\Vert u\Vert_{\cal H}^2 +\Vert Au\Vert_{\cal H}^2}$.  
\begin{definition}
The  complex number $z$  is in resolvent  set $\rho(A)$ of $A$  if and only if $A-z\1$ is invertible   from $D(A)$ into ${\cal H}$ and  $(A-z\1)^{-1}\in {\cal L}({\cal H})$  where  ${\cal L}({\cal H}$ is the Banach space  of linear and continuous maps  in ${\cal H}$.\\
The spectrum $\sigma(A)$ is defined as  $\sigma(A) = \C\backslash\rho(A)$
\end{definition}
\begin{definition}
Fix $\varepsilon >0$. The $\varepsilon$-spectrum $\sigma_\varepsilon(A)$ of $A$ is defined as follows.
A complex number $z\in\sigma_\varepsilon(A)$  if and only if $z\in\sigma(A)$ or if 
$\Vert(A-z\1)^{-1}\Vert_{{\cal L}({\cal H})} > \varepsilon^{-1}$.\\
It is  convenient   to write  $\Vert(A-z\1)^{-1}\Vert_{{\cal L}({\cal H})}=\infty$  if $z\in\sigma(A)$
 and denote $A-z = A-z\1$.
\label{Def1}
\end{definition}
There  are  several  equivalent  definitions of $\sigma_\varepsilon(A)$  for details see  the introduction of the book \cite{trem}. The following characterization is useful  for numerical  computations.  \\
 Assume that  dim${\cal H} <+\infty$.  Recall that the  singular values  for $A\in{\cal L}({\cal H})$    are  the eigenvalues
  of  the non negative matrix $\sqrt{A^*A}:=\vert A\vert$.   Denote $s(A) = \sigma(\vert A\vert)$.
  
  \begin{proposition}
For any matrix $A$ we have $z\in\sigma_\varepsilon(A)$  if and only if  $s_{\rm min}(A-z)]<\varepsilon$, where we have denoted
 $s_{\rm min}(A) := \min[s(A)] $.
\end{proposition}
{\bf Proof}  It is known that $\Vert A\Vert = s_{\rm max}(A)$ for every $A\in{\cal L}({\cal H})$.  But $AA^*$ and $A^*A$ have the same  non zero eigenvalues, so   if $A$ is invertible we have we have  $\Vert A^{-1}\Vert = \frac{1}{s_{min}(A)}$ and the proposition follows.
$\square$.

\subsection{Pseudospectra for quadratic pencils}
Our numerical  computations (see hereafter Section 5) show that the spectra  of quadratic   pencils is much more  unstable than the spectra of linear pencils (rotated harmonic oscillator, see \cite{Davies}).\\
Let us  recall  the basic  definitions and properties  concerning    pseudospectra for  quadratic pencils.
A more  general setting is  explained in \cite{hiti1, hiti2} for  pencils  of matrices.

The following result gives  an idea about the pseudospectra of the Schr\"odinger  pencil 
$L(\lambda) = -\frac{d^2}{dx^2} + (x^2-\lambda)^2$ :
\begin{theorem} \cite{chr3}
Assume that $\theta\in\R$,  $0<\vert\theta\vert \leq\frac{\pi}{2}$ and denote $\lambda_0= \rho e^{i\theta}$.  Then there exists $C<+\infty$, $\delta >0$ and for every  $\rho \geq 1$ a Schwartz  function $g$, $\|g\|=1$ such  that :
 \begin{equation}
  \| L(\lambda_0)g\|_{\mathcal{L}^2(\mathbb{R})}\leq C e^{-\delta \rho^{\frac{3}{2}}}
 \end{equation}
\end{theorem}
  i.e. for $\rho$ large enough the complex number $\lambda_0$ is, in some sense,  an almost eigenvalue   or     a pseudospectral  point  of $L(\lambda)$.   On the line  of     direction $\theta\in ]0, \pi/2]$ we have 
   for $\vert\lambda\vert$ large enough, 
  $$
  \Vert L(\lambda)^{-1}\Vert \geq \frac{1}{C}{\rm e}^{\vert\lambda\vert^{3/2}}
  $$
   In order  to   capture   more details  for the   localization  in the  complex plane  of large   modulus pseudospectral points   of   $L(\lambda)$  we  can consider   the  following  tentative definition  of    pseudospectra.\\
  Let  us consider  a quadratic pencil $L(\lambda)$   satisfying  the  assumptions of Theorem \ref{gene}.
  \begin{definition}
  Let $\varepsilon >0$, $\delta\geq 0$, $\mu>3/2$.  Define  the  pseudospectra of order $(\varepsilon, \delta, \mu)$   as  follows
 \begin{equation}
  {\rm Sp}_{\varepsilon,\delta,\mu}[L] = \{\lambda\in\C,\; \Vert L(\lambda)^{-1}\Vert \geq\varepsilon^{-1}\exp(\delta\vert\lambda\vert^\mu)\}
\label{Pseudo1}
 \end{equation}
  \label{Def2}
  \end{definition}
For  $\delta =0$  we recover  the  definition  given  by  Threfeten.\\
It is clear  that $\lambda\in  
{\rm Sp}_{\varepsilon,\delta,\mu}[L]$  if and only if  there exists $u\in D(L_0)$, $u\neq 0$,   such  that
  $$
\Vert L(\lambda)u\Vert \leq \varepsilon\exp(-\delta\vert\lambda\vert^\mu)\Vert u\Vert
  $$

\begin{remark}
Later we shall compute pseudospectra with this definition and see how it behaves according the parameter $0 \le \mu < \infty$.
\end{remark}

\section{Eigenvalues computation with spectral methods and finite difference discretization}

The aim of this section is to present the numerical computation of the spectrum of linear operator with quadratic dependence (quadratic pencil), see \Ref{Deq2}~:
$$
L(\lambda) = L_0 + \lambda L_1 + \lambda^2 
$$
where $L_0$ and $L_1$ are operators on some Hilbert space ${\mathcal H}$. 
So we are interested to solve the following nonlinear eigenvalue problem~:
$$
L(\lambda) u = 0 \ \ , \ \ \lambda\in \BigC,\;\;u\in {\mathcal H}.  
$$
In a first step, in order to validate the numerical approaches proposed, we consider the rotated harmonic oscillator in $L^2(\R)$(see Davies \cite{Davies})~:
\begin{equation}
- h \Frac{d^2}{dx^2} + c x^2 
\label{FJeq2}
\end{equation}
where $h$ is a real positive parameter and $c$ is a complex number with positive real and imaginary parts, $c=\exp(i\alpha)$, for $0 \le \alpha < \pi/2$.\\

Here after, for each operator considered, we compute spectra and pseudospectra and we discuss the numerical results obtained. 

\subsection{The rotated harmonic oscillator}

\subsubsection{Eigenvalue computations with Hermite spectral method (unbounded domain)}

We look for an approximation of $\lambda \in \BigC$ solution of the following linear eigenvalue problem ($h=1$)~:
\begin{equation}
-\Frac{d^2u}{d x^2}(x) + c x^2 u(x) = \lambda u(x) \ \ , \ \ x \in \BigR 
\label{eq3ori}
\end{equation}
Here the computational domain is unbounded ($\Omega = \BigR$). So we use a spectral Galerkin method using Hermite functions (see Appendix A) {\it i.e.} we look for an approximation~:
\begin{equation}
u_N(x) = \displaystyle \sum_{k=0}^N \tilde{u}_k \varphi_k (x)
\label{eq3bis}
\end{equation}
of $u$ such that~: 
$$
\la -\Frac{d^2u_N}{dx^2} + c x^2 u_N - \lambda_N u_N,\varphi_{l} \ra  = 0 \ \ , \ \ l=0,\ldots,N
$$
where $\la, \ra$ is the scalar product in $L^2(\R)$ (method of weighted residuals, MWR, see for example, \cite{GO}, \cite{Cal}). Using the orthogonality properties \Ref{crerec} of the Hermite functions in $L^2(\R)$ and the relations (\ref{psphix}) we obtain the following eigenvalue problem~:
\begin{equation}
{\mathcal A}_N U_N = \lambda_N U_N
\label{eq3ter}
\end{equation}
with ${\mathcal A}_N$ the square tridiagonal symmetric matrix of order $N+1$ such that ${\mathcal A}_N(k,k-2) = (c-1)\sqrt{k(k-1)}$, ${\mathcal A}_N(k,k) = (c+1)(2k+1)$ and $U_N$ is the vector containing the coefficients $\tilde{u}_k$, $k=0,\ldots N$ of $u_N$. \\

For the numerical computation of the spectrum of ${\mathcal A}_N$ we use the function ZGEEV of the library LAPack. \\

We recall that for the continuous operator \Ref{eq3ori} the eigenvalues are (see \cite{Zworski})~:
\begin{equation}
\exp (i\alpha/2) (2n+1) \ \ , \ \ n = 0,1,\ldots
\label{eq3qua}
\end{equation}
inducing that the eigenvalues, in the complex plane, are aligned on a straight with a slope $\Frac{\lambda_i}{\lambda_r} = \tan(\alpha/2)$, where $\lambda_r$ (resp. $\lambda_i$) is the real (resp. imaginary) part of the eigenvalues $\lambda$. Here we have chosen $\alpha = \pi/4$ so $\Frac{\lambda_i}{\lambda_r} = \tan(\pi/8) = 0.4$.\\
 
\begin{figure}
\centering
\vspace{4cm}
\includegraphics[width=4in]{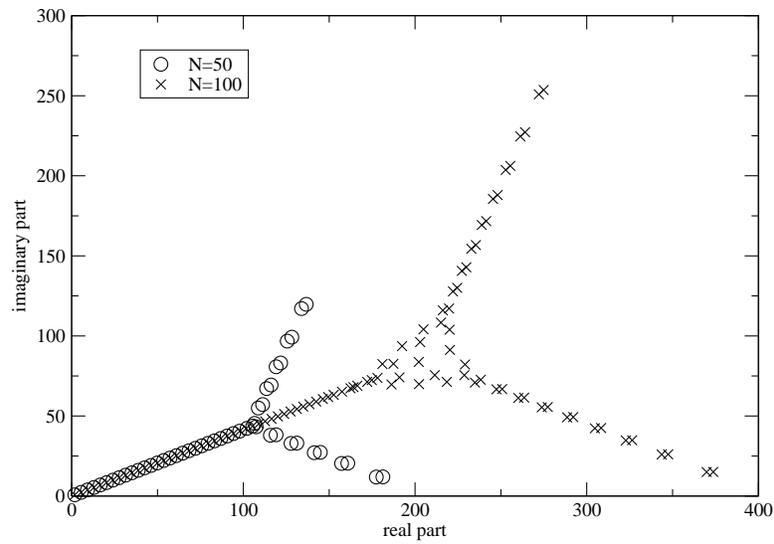}
\caption{ \label{Figure1} Spectrum of the matrix ${\mathcal A}_N$ associated with the eigenvalue problem \Ref{eq3ter} (Hermite spectral method) for $N=50$ and $N=100$, $c=\exp(i\alpha)$ with $\alpha=\pi/4$. }
\end{figure}

Now we present the numerical results obtained with the Hermite spectral method. On Figure \ref{Figure1} we can see the spectrum of the matrix ${\mathcal A}_N$ associated with the eigenvalue problem \Ref{eq3ter} for $N=50$, $c=\exp(i\alpha)$ with $\alpha=\pi/4$. We can see that the slope $\lambda_{N,i}/\lambda_{N,r} = \tan (\alpha/2)$ is obtained for $\vert \lambda_{N,r} \vert \le 100$. Then a bifurcation appears in the spectrum, which is in agreement with \cite{Davies}, \cite{Zworski}. If we choose a larger 
value of $N$, for example $N=100$, the slope of $\tan(\alpha/2)$ appears for larger value of $\vert \lambda_{N,r}\vert \le 200$ (see Figure \ref{Figure1}) which is in agreement with the fact that $u_N$ converges to $u$ when $N$ increases (see \Ref{err1}, \Ref{err2}). 

\subsubsection{Eigenvalue computations with finite difference method (bounded domain)}

The rotated harmonic oscillator is defined for functions $u\in H^2(\R)$ such that $x^2 u(x) \in L_2(\BigR)$. So $u(x)$ decreases when $x^2$ increases and we want to consider the problem \Ref{eq3ori} on a bounded domain, with homogeneous Dirichlet boundary conditions~: 
\begin{equation}
\left\{
\begin{array}{l}
-\Frac{d^2u(x)}{dx^2} + c x^2 u(x) = \lambda u(x) \ \ , \ \ x \in \Omega \\
u(\pm L) = 0 
\end{array}
\right.
\label{eq4}
\end{equation}
where $\Omega = (-L,+L)$, $L$ being chosen sufficiently large. More precisely, if we retain $N$ modes in the Hermite development \Ref{eq3bis}, the Hermite function of highest degree is $\varphi_N$ and the zeroes $h_n$ of $\varphi_N $
verify (see \cite{ADGR})~:
$$
h_n \le \sqrt{2N-2} \ \ , \ \ n = 1, \ldots N
$$
So the size of the containment area is $2L = 2\sqrt{2N-2}$ and we retain as value for the bounded domain $\Omega$~:
\begin{equation}
L = \sqrt{2N-2}
\label{eq6}
\end{equation}
To obtain a numerical approximation $\lambda_N$ of the eigenvalues of problem \Ref{eq4} we discretize the second order spatial derivative using a standard second order centered finite difference scheme~:
$$
\Frac{d^2u (x_j)}{dx^2} = \Frac{u(x_{j+1})-2u(x_j)+u(x_{j-1})}{\Delta x ^2} + o(\Delta x^2)
$$
where $\Delta x = 2L/N$ is the spatial step of the meshgrid $x_j = -L + j\Delta x$, $j=0,\ldots N$, on the domain $\Omega$. So we obtain the following linear eigenvalue problem to solve~:
\begin{equation}
{\mathcal A}_N U_N = \lambda_N U_N
\label{eq7}
\end{equation}
where ${\mathcal A}_N$ is a tridiagonal symmetric matrix of order $N-1$ such that ${\mathcal A}_N(k,k-1) = {\mathcal A}_N (k,k+1) = -\Frac{1}{\Delta x^2}$, ${\mathcal A}_N(k,k) = \Frac{2}{\Delta x^2} +cx_k^2$ and $U_N$ is the vector containing the approximations $u_N(x_j)$ of $u(x_j)$, $j=1,\ldots N-1$ ($u_N(x_0) = u_N(x_N) = 0$).\\

As previously for the Hermite spectral method, we use the function ZGEEV of the library LAPack for the numerical computation of the spectrum. \\

\begin{figure}
  \centering
\includegraphics[width=4in]{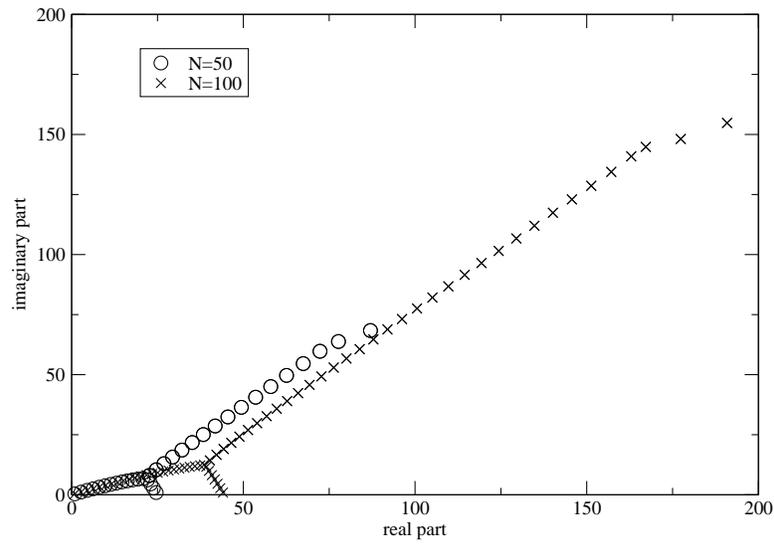}
\caption{ \label{Figure2} Spectrum of the matrix ${\mathcal A}_N$ of the eigenvalue problem \Ref{eq7} (finite difference scheme) obtained for $N=50$ ($L=10$), $N=100$ ($L=15$) and $c=\exp(i\alpha)$ with $\alpha=\pi/4$.}
\end{figure}

Now we present the numerical results obtained with the method based on finite difference discretization. As for the Hermite spectral method, we have chosen $N=50$ and $N=100$. 
On Figure \ref{Figure2} we present the spectrum of the matrix ${\mathcal A}_N$ of the eigenvalue problem \Ref{eq7} obtained for $N=50$, $L=10$ in accordance with \Ref{eq6} and $c=\exp(i\alpha)$ with $\alpha=\pi/4$. We can see that the slope $\lambda_{N,i}/\lambda_{N,r} = \tan (\alpha/2)$ is obtained for $\vert \lambda_{N,r} \vert \le 25$. Then, as for the Hermite spectral method, a bifurcation appears in the spectrum. If we choose a larger 
value of $N$, for example $N=100$ and $L=15$ following \Ref{eq6}, the slope of $1/2$ appears in the spectrum for larger value of $\vert \lambda_{N,r} \vert \le 75$ (see Figure \ref{Figure2}) which is in agreement with the fact that accuracy of the difference scheme increases with $N$. \\

\begin{figure}
  \centering
\includegraphics[width=4in]{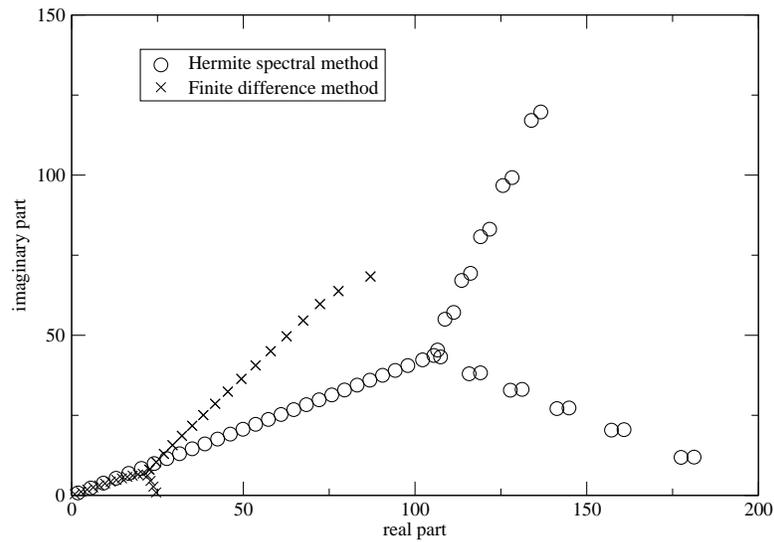}
\caption{ \label{Figure3} For $N=50$, comparison of the numerical results obtained with the Hermite spectral method and with the finite difference scheme.}
\end{figure}

Now if we compare, for a same value of $N$ ($N=50$) the numerical results obtained with Hermite spectral method and with finite difference scheme, we can see on Figure \ref{Figure3} that the slope $\Frac{\lambda_{N,i}}{\lambda_{N,r}} = \tan(\alpha/2)$, which is in agreement with the continuous operator (see \Ref{eq3qua}), appears for larger values of $\lambda_{N,r}$ with the spectral method than with the finite difference method. This is coherent with the fact that the Hermite spectral method is more accurate than the finite difference method for a same value of the parameter $N$ (spectral accuracy due to the fast decrease, in modulus, of the coefficients $\tilde{u}_k$ when $k$ increases, see Proposition \Ref{he}). \\


Now we try to analyze the bifurcation phenomenon appearing on the spectrum for eigenvalues with large real part $\lambda_{N,r}$. When we discretize with a finite difference scheme, we consider that $x$ is constant over one spatial step $\Delta x = \Frac{2L}{N}$. So, in a first step, we consider an operator deduced from the rotated harmonic oscillator in which $x^2$ is chosen constant equal to $b^2$ over all the domain $\Omega$. So we obtain the following problem deduced from \Ref{eq4}~:
\begin{equation}
-\Frac{d^2u}{d x^2} + c b^2 u = \lambda u \ \ , \ \ x \in \Omega \\
\label{eq8}
\end{equation}
If we consider periodic boundary conditions $u(-L)=u(+L)$, we look for eigenfunctions of \Ref{eq8} such as~:
\begin{equation}
u(x) = \hat{u}_k \exp(ik'x)
\label{eq9}
\end{equation}
with the wavenumber $k' = \Frac{k\pi}{L}$, $k=0,\ldots N-1$. Substituting \Ref{eq9} in \Ref{eq8} we obtain~:
$$
\lambda = k'^2 + c b^2 
$$
So, since $c=\exp(i\alpha)$ we have~:
$$
\left\{
\begin{array}{l}
\lambda_r = k'^2 + \cos(\alpha) b^2 \\
\lambda_i = \sin(\alpha) b^2 
\end{array}
\right.
$$
We can see that $\lambda_i$ is constant and that $\lambda_r$ depends of the wavenumber $k'$. \\

\begin{figure}
  \centering
  \includegraphics[width=4in]{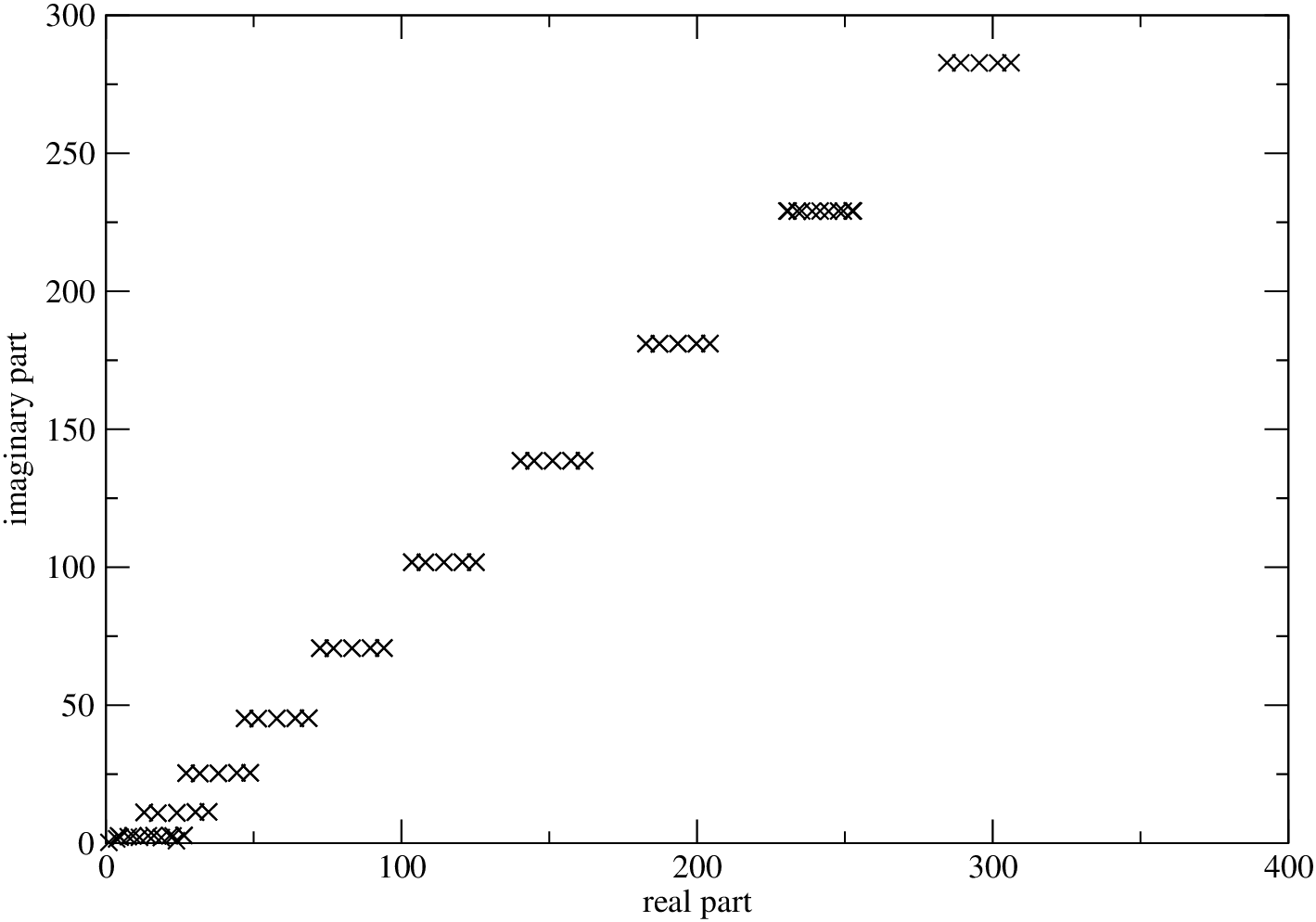}
\caption{ \label{Figure4a} Eigenvalues computed with the finite difference scheme, corresponding to $N=100$, $L=20$, $\alpha = \pi / 4$ and $N_b=5$. }
\end{figure}

\begin{figure}
  \centering
 \includegraphics[width=4in]{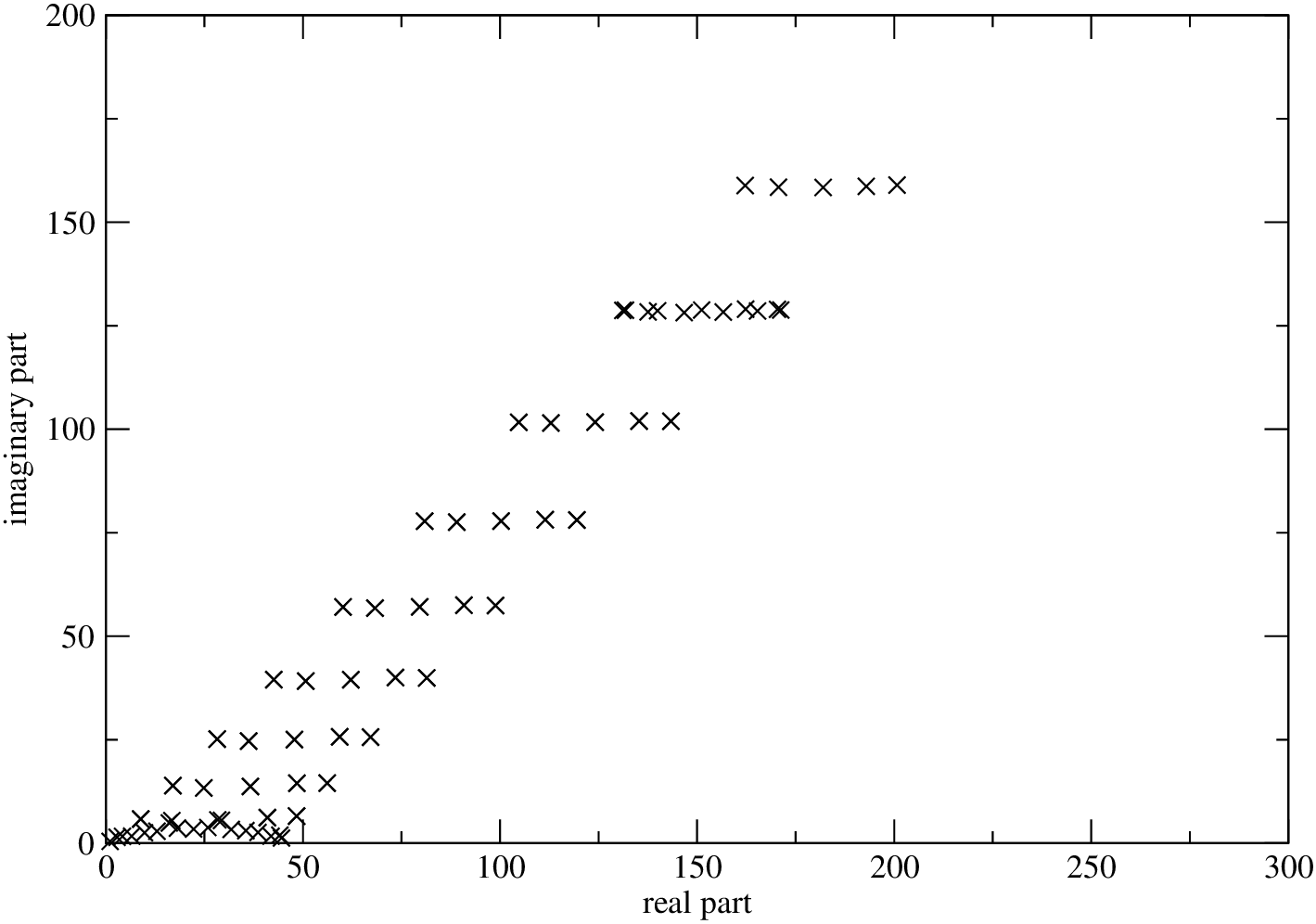}
\caption{ \label{Figure4b} Eigenvalues computed with the finite difference scheme, corresponding to $N=100$, $L=15$, $\alpha = \pi / 4$ and $N_b=5$.}
\end{figure}

Now we consider that $x$ is constant over $N_b$ spatial steps $\Delta x$, so in the rotated harmonic oscillator we replace $x^2$ with $b(x)^2$ where $b(x) = b_j = -L + j N_b \Delta x$ for $x \in [-L+jN_b\Delta x, -L+(j+1)N_b\Delta x[$, $j = 0,\ldots N/N_b-1$. We obtain~:
\begin{equation}
\left\{
\begin{array}{l}
\lambda_r= k'^2 + \cos(\alpha) b_j^2 \\
\lambda_i = \sin(\alpha) b_j^2 
\end{array}
\right.
\label{eq11}
\end{equation}
with the wavenumber $k' = \frac{k\pi}{L}$, $k=0,\ldots N_b-1$.
So the spectrum is constituted of different steps, each step corresponding to $N_b$ eigenvalues $\lambda$, with $\lambda_i$ constant while $\lambda_r$ is wavenumber dependent. 
We can observe on \Ref{eq11} that for $k'=0$ the corresponding eigenvalues $\lambda$, $j=0,\ldots,N/N_b-1$ are aligned on a straight with a slope $\lambda_i/\lambda_r = \tan(\alpha)$. This can be seen on Figure \ref{Figure4a}, corresponding to $N=100$, $L=20$, $\alpha = \pi / 4$ and $N_b=5$. On Figure \ref{Figure4b}, corresponding to $N=100$, $L=15$, $\alpha = \pi / 4$ and $N_b=5$ ($\Delta x$ is decreased in comparison with Figure \ref{Figure4a}), we can see that some numerical artefacts appear on the computation of the eigenvalues $\lambda$ having small modulus. 
In order to avoid this, we try to impose in the spectrum that two consecutive steps, corresponding to two different values of $b_j$, are not recovered for the real part $\lambda_r$. So we must have~:
$$
(\Frac{N}{2L})^2 \le \cos(\alpha) N_b^2 (\Delta x )^2
$$
So, since $\Delta x = \Frac{2L}{N}$ we deduce the following inequality~: 
\begin{equation}
(\Frac{N}{L})^2 \le 4 N_b \sqrt{\cos(\alpha)}
\label{eq12}
\end{equation}
which is a constraint on $\Delta x^{-1}$. \\

\begin{figure}
  \centering
 \includegraphics[width=4in]{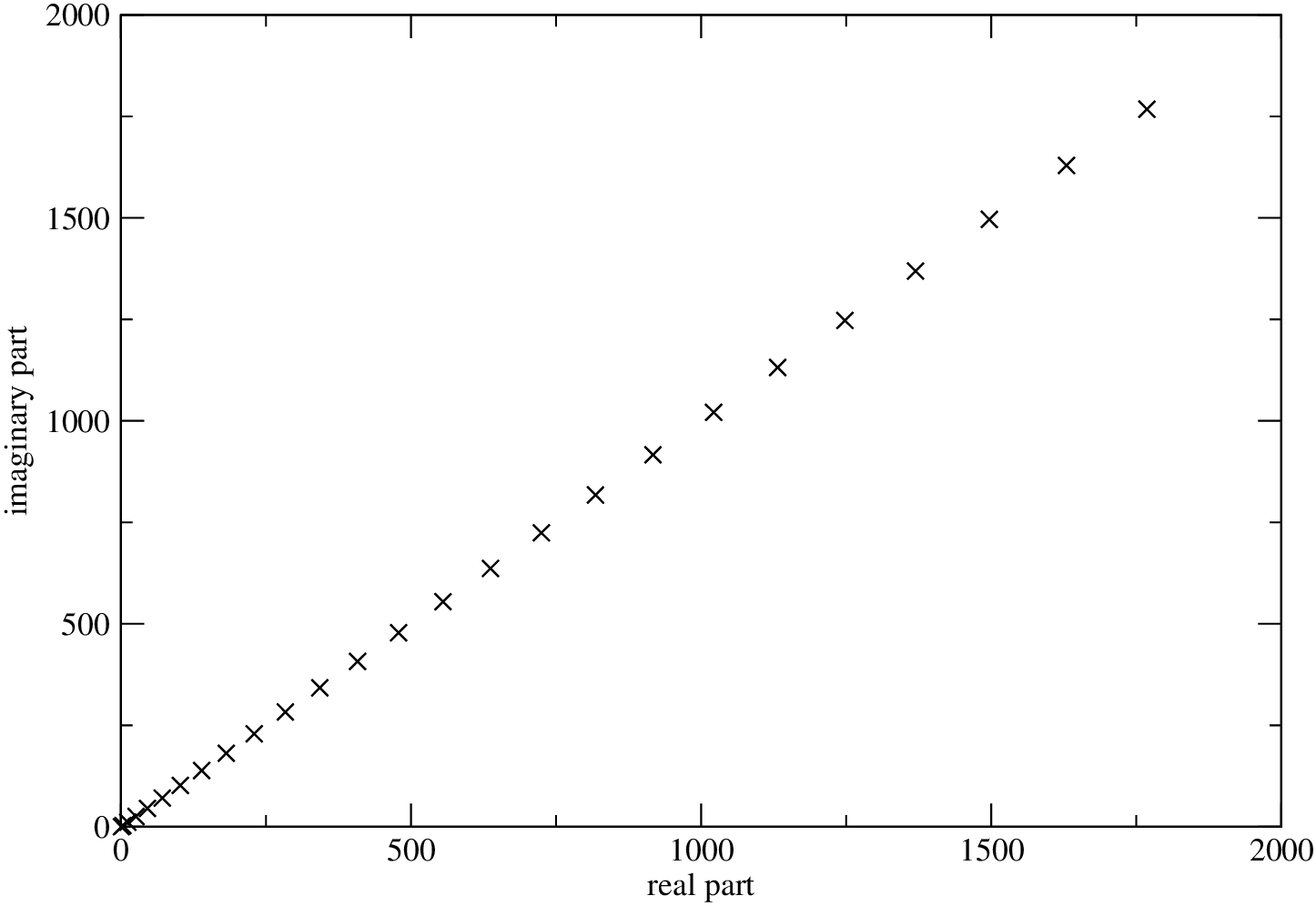}
\caption{ \label{Figure5a} Eigenvalues computed with the finite difference scheme \Ref{eq7} for $\alpha = \pi/4$, $N=50$ and $L=50$.}
\end{figure}

\begin{figure}
  \centering
 \includegraphics[width=4in]{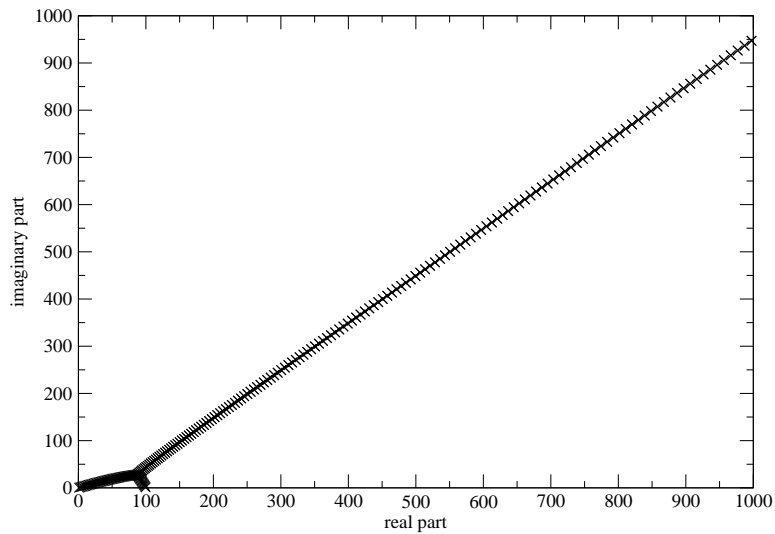}
\caption{ \label{Figure5b} Eigenvalues computed with the finite difference scheme \Ref{eq7} for $\alpha = \pi/4$, $N=500$ and $L=50$.} 
\end{figure}

\begin{figure}
  \centering
 \includegraphics[width=4in]{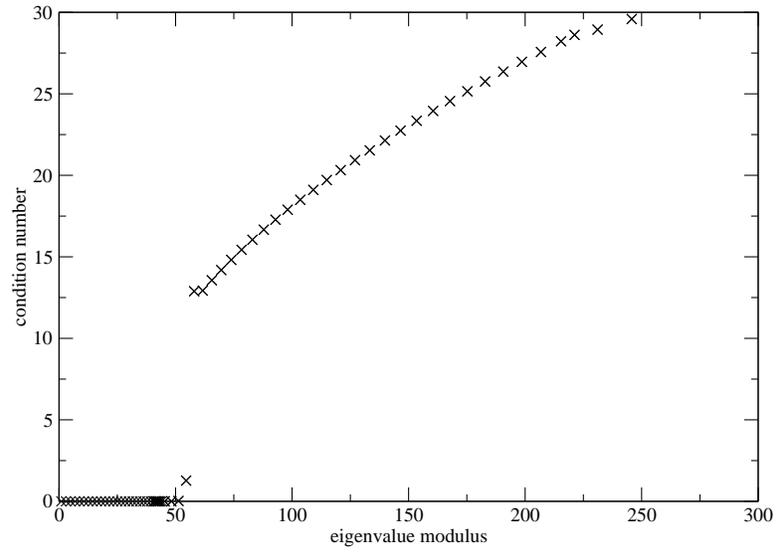}
\caption{ \label{Figure6} Condition number of the eigenvalues $\lambda_N$ in function of the modulus of the eigenvalues, for $N=100$ and $L=15$.}
\end{figure}

As it has been said previously for the finite difference scheme we have $N_b = 1$. So, in agreement with \Ref{eq11} we expect that the eigenvalues computed with the finite difference scheme \Ref{eq7} are aligned on a straight with a slope $\lambda_i/\lambda_r = \tan(\alpha) = 1$ for $\alpha = \pi/4$. This is what we obtain if we choose $N=50$ and $L=50$ (see Figure \ref{Figure5a}). We can note that with this choice of the parameters, the inequality \Ref{eq12} is satisfied. Now, in order to test the convergence of the finite difference scheme we reduce the spatial step $\Delta x$. So we choose $N=500$ and $L=50$ (see Figure \ref{Figure5b}). With this choice of the parameters, the inequality \Ref{eq12} is not satisfied. We can see that some numerical artefacts appear near the origin, where we can observe a slope $\Frac{\lambda_{N,i}}{\lambda_{N,r}} = \tan(\alpha/2)$, which is in agreement with \Ref{eq3qua}. This can be interpreted as an intermediate slope between the slope $\Frac{\lambda_{N,i}}{\lambda_{N,r}} = \tan(\alpha)$ and the slope $\Frac{\lambda_{N,i}}{\lambda_{N,r}} = \tan (0) = 0$ of each step. 
\begin{remark}
Let $H(b_j)$ the operator deduced from the rotated harmonic oscillator, $H(b_j) = -\Frac{d^2}{dx^2} + c b_j^2$ and $\lambda_{j,i}$, $u_{j,i}$ the eigenvalues and eigenvectors associated~: $H(b_j) u_{j,i} = \lambda_{j,i} u_{j,i}$. 
We denote $u = \di \sum_{j=0}^{N/N_b-1} u_j \BigI_{\left[x_j,x_{j+1}\right]}$ and $b = \di \sum_{j=0}^{N/N_b-1} b_j \BigI_{\left[x_j,x_{j+1}\right]}$, with $x_j = -L+j N_b \Delta x$ and $\BigI_{\left[x_j,x_{j+1}\right]}$ the characteristic function associated with the interval $\left[x_j,x_{j+1}\right]$. We consider the operator $H(b) = \oplus H(b_j)$. If $\lambda$ is an eigenvalue of $H(b)$, so there exists $(j,i)$ such that $\lambda = \lambda_{j,i}$.  
\end{remark}
Now, in order to study the numerical instability of the finite difference scheme in function of the meshgrid $x_j$, $j=0,\ldots N$, we consider a small perturbation on each point of the grid, $x_j + \varepsilon$, $j=0,\ldots N$, where $\varepsilon$ is a small parameter. The matrix ${\mathcal A}_N$ (see \Ref{eq7}) is replaced with the matrix~:
$$
{\mathcal A}_{N,\varepsilon} = {\mathcal A}_N + \varepsilon {\mathcal E}_N
$$
where ${\mathcal E}_N$ is the diagonal matrix of order $N-1$ such that ${\mathcal E}_N (k,k) = 2 \exp(i\alpha) x_k$ (we have neglected the terms in $\varepsilon^2$). If we compare the eigenvalues $\lambda_N$ of the matrix ${\mathcal A}_N$ with the eigenvalues $\lambda_{N,\varepsilon}$ of ${\mathcal A}_{N,\varepsilon}$ we have~:
$$ 
\mathcal{A}_{N,\varepsilon} U_{N,\varepsilon} = ({\mathcal A}_N + \varepsilon {\mathcal E}_N) U_{N,\varepsilon} = \lambda_{N,\varepsilon} U_{N,\varepsilon}
$$
where $U_{N,\varepsilon}$ is a right eigenvector of $A_{N,\varepsilon}$. So we deduce (see \cite{Quarteronietal})~:
$$
{\mathcal A}_N \Frac{d U_{N,\varepsilon}}{d \varepsilon} (\varepsilon) + {\mathcal E}_N U_{N,\varepsilon} (\varepsilon) + \varepsilon {\mathcal E}_N \Frac{d U_{N,\varepsilon}}{d \varepsilon} (\varepsilon) = \Frac{d \lambda_{N,\varepsilon}}{d\varepsilon} (\varepsilon) U_{N,\varepsilon} (\varepsilon) + \lambda_{N,\varepsilon} (\varepsilon) \Frac{d U_{N,\varepsilon}}{d \varepsilon} (\varepsilon)
$$
For $\varepsilon = 0$ we obtain~:
$$
{\mathcal A}_N \Frac{d U_{N,\varepsilon}}{d \varepsilon} (0) + {\mathcal E}_N U_{N,\varepsilon} (0) = \Frac{d \lambda_{N,\varepsilon}}{d\varepsilon} (0) U_{N,\varepsilon} (0) + \lambda_{N,\varepsilon} (0) \Frac{d U_{N,\varepsilon}}{d \varepsilon} (0)
$$
If we multiply on the left the previous equality with $V_N$ a left eigenvector of ${\mathcal A}_N$ we obtain~:
\begin{equation}
\Frac{d \lambda_{N,\varepsilon}}{d \varepsilon}(0) = \Frac{V_N^{\star} {\mathcal E}_N U_N}{V_N^{\star} U_N}
\label{eq14}
\end{equation}
where $V_N^{\star} = {\overline V_N}^t$ and $U_N$ is a right eigenvector of ${\mathcal A}_N$. The equality \Ref{eq14} measures the sensivity of the eigenvalue $\lambda_N$ of the matrix ${\mathcal A}_N$ in function of a perturbation $\varepsilon$ on the meshgrid (condition number of the eigenvalue $\lambda_N$). On Figure \ref{Figure6} we have represented the condition number of the eigenvalues $\lambda_N$ in function of the modulus of the eigenvalues, for $N=100$ and $L=15$. We can see that the condition number is small for eigenvalues with small modulus and then it increases with the modulus. However, the values are small in comparison with the results obtain for a nonlinear eigenvalue problem (see Section 5.2.2, Figure \ref{Figure11}), which implies that the numerical computation is stable if we consider a perturbation on the meshgrid points. 

\subsubsection{Pseudospectra}

In this subsection we present  numerical pseudospectra computations  for the rotated harmonic oscillator.  Notice that  a   theoretical  analysis  of this problem  has been performed in \cite{ps}.   In a first step we consider the matrix from the Hermite spectral method \Ref{eq3ter} and in a second step the matrix from the finite difference scheme \Ref{eq7}. 
It is known that the numerical computation of the pseudospectra is more stable than for the spectra (see Section 4). \\

To obtain the pseudospectra, following Definition \ref{Def1} we look for $z \in \BigC$ such that $\vert \vert ( {\mathcal A}_N -z I_N )^{-1} \vert\vert = s^{-1}_{\rm min} ({\mathcal A}_N-z \1_N)$ is large, {\it i.e.} the distance of $z$ to the spectrum of ${\mathcal A}_N$ is small~:
\begin{equation}
s_{\rm min} ({\mathcal A}_N-z \1_N)\le \varepsilon
\label{eq15}
\end{equation}
where $\vert\vert - \vert\vert$ is the matricial norm associated with the Euclidean norm, $\varepsilon$ is a small parameter, $\1_N$ is the identity matrix and $s_{\rm min} ({\mathcal A}_N-z \1_N)$ is the smallest singular value of the matrix ${\mathcal A}_N-z \1_N$. So we consider a mesh on the complex plane. For each point $z$ of the mesh we compute the singular value of ${\mathcal A}_N-z \1_N$, using the function ZGESVD of the LAPack Library. \\

\begin{figure}
  \centering
\includegraphics[width=4in]{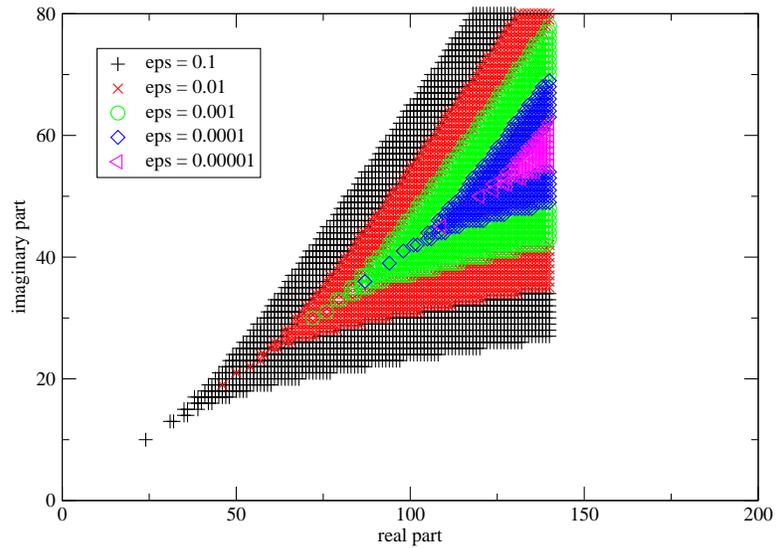}
\caption{ \label{Figure7a} Computation of the pseudospectra \Ref{eq15} of the matrix ${\mathcal A}_N$ corresponding to the Hermite spectral method for $\alpha = \pi/4$, $N=100$. }
\end{figure}

\begin{figure}
  \centering
 \includegraphics[width=4in]{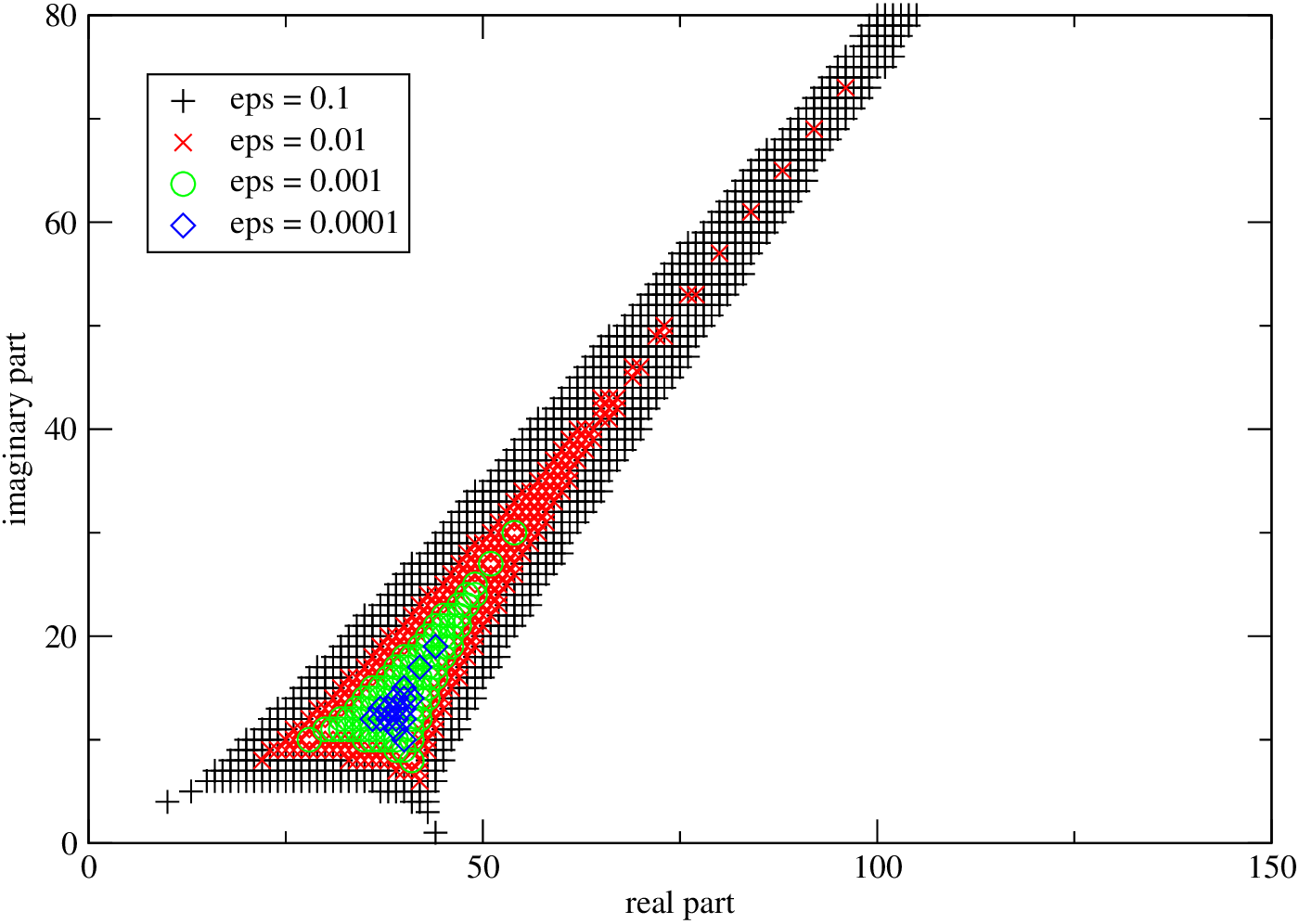}
\caption{ \label{Figure7b} Computation of the pseudospectra \Ref{eq15} of the matrix ${\mathcal A}_N$ corresponding to the finite difference scheme for $\alpha = \pi/4$, $N=100$ and $L=15$. }
\end{figure}

For the computation of the pseudospectra \Ref{eq15}, we have retained complex values $z$ lying on the meshgrid in the area of the complex plane corresponding to $[0,140] \times [0,80]$. The step retained is $dx=1$ and $dy=1$ in the real and imaginary directions. On Figure \ref{Figure7a} (resp. Figure \ref{Figure7b}) we can see the computation corresponding to the matrix ${\mathcal A}_N - z \1_N$, with ${\mathcal A}_N$ corresponding to the Hermite spectral method (resp. finite difference scheme). The choice of the parameters are $\alpha = \pi/4$, $N=100$ for unbounded and bounded domains, $L=15$ for the bounded domain. We can see on these two figures that the spectrum of the continuous operator (slope equal to $\lambda_i/\lambda_r = \tan (\alpha/2)= 0.4$ (see \Ref{eq3qua}) is contained in the area of the pseudospectra corresponding to the smallest values of the parameter $\varepsilon$, {\it i.e.} in the area where the distance of $z$ to the eigenvalues of matrix ${\mathcal A}_N$ is the smallest. This is especially true for the Hermite spectral method. \\

Now we consider here the computation of the pseudospectra based on Definition \ref{Def2} (see \Ref{Pseudo1}) instead of Definition \ref{Def1} as 
previously. So we look for $z \in \BigC$ such that~:
\begin{equation}
\vert \vert ({\mathcal A}_N -z \1_N )^{-1}\vert\vert =  s^{-1}_{\rm min} ({\mathcal A}_N-z \1_N)\ge \varepsilon^{-1} \exp(\delta \vert z \vert ^\mu)
\label{Pseudo2}
\end{equation}
where $\varepsilon$, $\delta$ and $\mu$ are real parameters. For $\delta = 0$ we retrieve \Ref{eq15}. The majoration in \Ref{Pseudo2} depends on $\vert z \vert$, in opposition with \Ref{eq15}~: the constraint appearing in \Ref{Pseudo2} is stronger when $\vert z \vert$ is increased.  
On Figure \ref{Figure100} we present the pseudospectra corresponding to the matrix ${\mathcal A}_N - z \1_N$, where ${\mathcal A}_N$ is the matrix obtained with the finite difference scheme, for $\delta = 0.5$, $\mu = 0.5$ and for different values of the parameter $\varepsilon$. As before we have retained complex values $z$ lying on the meshgrid in the area of the complex plane corresponding to $[0,140] \times [0,80]$. The step retained is $dx=1$ and $dy=1$ in the real and imaginary directions. The numerical results obtained are in agreement with the results presented on Figure \ref{Figure7b}.

\subsection{Nonlinear eigenvalue problems}

In this section we consider the following operator~:
\begin{equation}
L_{a} (\lambda) = - \Frac{d^2}{dx^2} + x^4 - 2 a \lambda x^2 + \lambda^2
\label{FJeq16}
\end{equation}
where $a$ is a real parameter. We want to solve the following problem~:
\begin{equation}
L_{a} (\lambda) u = 0 \ \ , \ \ x \in \BigR
\label{FJeq17}
\end{equation}
For $a=1$ we recover the problem \Ref{Deq3}. \\

The problem \Ref{FJeq17} car be reformulated as an eigenvalue problem. Indeed, if we set $v = \lambda u$ we can rewrite \Ref{FJeq17} as~:
\begin{equation}
{\mathcal  A}_{a} \left(
\begin{array}{l} 
u \\
v
\end{array}
\right) 
= \lambda 
\left(
\begin{array}{l}
u \\
v
\end{array}
\right) 
\label{FJeq22bis}
\end{equation}
where~:
$$
                      {\mathcal  A}_{a} = \left(
                            \begin{array}{cc}
                            0 & \1\\
                            -L_0&-L_{1_a}\\
                            \end{array}
                            \right)
$$
with the operators $L_0 = - \Frac{d^2}{dx^2} + x^4$ and $L_{1_a} = - 2 a x^2$. 

\subsubsection{Eigenvalue computations with Hermite spectral method (unbounded domain)}

We look for an approximation $u_N$ of $u$ such that $L_a(\lambda) u_N = 0$, with $u_N = \displaystyle \sum_{k=0}^N \tilde{u}_k \varphi_k$, with $\varphi_k$ Hermite functions (spectral Galerkin approximation, see Appendix A). Then, in order to obtain $u_N$, we use a method of weighted residuals (MWR, see for example, \cite{GO}, \cite{Cal})~:
$$
(L_a(\lambda) u_N,\varphi_{l}) = 0 \ \ , \ \ l=0,\ldots,N
$$
where $(.,.)$  is the scalar product in $L^2(\R)$. Setting $v_N = \lambda u_N$, using the orthogonality properties of the Hermite function in $L^2(\R)$ and the relations (\ref{psphix}), we obtain the following eigenvalue problem~:
$$
{\mathcal A}_{a,N}
 \left(
\begin{array}{l} 
U_N \\
V_N
\end{array}
\right) 
= \lambda_N
\left(
\begin{array}{l}
U_N \\
V_N
\end{array}
\right) 
$$
which is an approximation of the eigenvalue problem \Ref{FJeq22bis}. $U_N$ (resp. $V_N$) is the vector containing the coefficients $\tilde{u}_k$ (resp. $\tilde{v}_k$) of $u$ (resp. $v$), $k=0,\ldots,N$. 
The matrix ${\mathcal A}_{a,N}$ is the square matrix of order $(2N+2)$~:
\beq
                      {\mathcal A}_{a,N} = \left(
                            \begin{array}{cc}
                            0 & \1_N\\
                            -L_{0_N}&-L_{1_{a,N}}\\
                            \end{array}
                            \right)
\label{FJeq17bis}
\eeq
where $L_{0_N} u_N = (L_0 u_N,\varphi_{l})$ and $L_{1_{a,N}} v_N = (L_{1_a} v_N,\varphi_{l})$, $l=0,\ldots,N$, with $L_0 u_N = -\Frac{d^2u_N}{dx^2} + x^4 u_N$ and $L_{1_a} v_N = -2 a  x^2 v_N$. \\


%
%
%
%
%
$L_{0_N}$ is a pentadiagonal symmetric matrix such that $L_{0_N}(j,j) = (\Frac{2j+1}{2}) + \Frac{1}{4}c_j$, $L_{0_N}(j,j-2) = - \Frac{1}{2}\sqrt{j(j-1)} + \Frac{1}{4}b_{j-2}$ and $L_{0_N}(j,j-4) = \Frac{1}{4} a_{j-4}$ for $j=0,\ldots,N$, where $a_j = \sqrt{j(j-1)(j-2)(j-3)}$, $b_j=(4j-2)\sqrt{j(j-1)}$ and $c_j=(6j^2+2j+3)$.\\

$L_{1_{a,N}}$ is a tridiagonal symmetric matrix such that $L_{1_{a,N}}(j,j) = - (2j+1)$ and $L_{1_{a,N}}(j,j-2) = - \sqrt{j(j-1)}$.\\

For the numerical computation of the spectrum of ${\mathcal A}_{a,N}$ we use the function DGEEV of the LAPack library.\\


For $a=1$, in order to analyze the spectrum of the continuous operator \Ref{Deq3}, we consider a simplified operator, deduced from the operator \Ref{FJeq16} for $a=1$, where $x$ is replaced with a real constant $b$. We obtain the following problem~:
\begin{equation}
\lambda^2 u - \Frac{d^2u}{dx^2} - 2 b^2 \lambda u + b^4 u = 0
\label{FJeq17ter}
\end{equation}
%
%
%
%
%
%
%

%
We look for a solution $u(x)$ of the problem \Ref{FJeq17ter} of the form $u = \tilde{u}_k \varphi_k$. Substituting in \Ref{FJeq17ter} and using the relations \Ref{psphix}, we obtain~:
\begin{equation}
\lambda_N^2 \varphi_k - \Frac{1}{2} \sqrt{k(k-1)} \varphi_{k-2} + (\Frac{2k+1}{2k+2}) \varphi_k - \Frac{1}{2} \sqrt{(k+1)(k+2)} \varphi_{k+2} - 2 b^2 \lambda_N \varphi_k + b^4 \varphi_k = 0
\label{FJeq18}
\end{equation}
Using the scalar product in $L^2(\R)$ of \Ref{FJeq18} with $\varphi_k$ we obtain~:
\begin{equation}
\lambda_N^2 + \Frac{2k+1}{2} - 2 b^2 \lambda_N + b^4 = 0
\label{FJeq19}
\end{equation}
We deduce from \Ref{FJeq19} that $\lambda_N = b^2 \pm i \sqrt{k+1/2}$, so~:
\begin{equation}
\left\{
\begin{array}{l}
\lambda_{N,r} = b^2 \\
\lambda_{N,i} = \pm \sqrt{k+1/2}
\end{array}
\right.
\label{FJeq20}
\end{equation}
The imaginary part $\lambda_{N,i}$ of $\lambda_N$ is wavenumber independent. From \Ref{FJeq20} it comes that the spectrum is contained in the part of the complex plane defined by $\lambda_{N,r} = b^2$ and $-\sqrt{N} \le \lambda_{N,i} \le \sqrt{N}$ since $k=0,\ldots, N$. \\

Now, on Figure \ref{Figure8} we present the spectrum of the matrix \Ref{FJeq17bis} for $N=50$ and $a=1$. Firstly we can note that, as for the rotated harmonic oscillator (see Figure \ref{Figure1}), a bifurcation appears in the spectrum when the modulus of the eigenvalues is increased (see also Figure \ref{Figure12}). Then theoretical results give that the eigenvalues of the continuous operator \Ref{FJeq16}, for $a=1$, are included in the two  sectors $\{\lambda\in\C,\; \vert\arg(\lambda)\vert \geq \Frac{\pi}{3}\}$ (see Section 3). We can see on Figure \ref{Figure8} that computed eigenvalues are not all included in these two sectors. This reflects numerical instabilities leading to spurious eigenvalues (spectral pollution, see \cite{DaviesPlum}). We can note that we have $-\sqrt{N} \le \lambda_{N,i} \le \sqrt{N}$, in agreement with the previous analyze when $x=b$ is constant (see \Ref{FJeq20}). Moreover, in the previous analyze we have $\lambda_{N,r} = b^2$. Here, for $N=50$ following \Ref{eq6} we deduce that the size of the containment domain is $2L$ with $L \simeq 10$ and, on Figure \ref{Figure8}, we can see that $0 \le \lambda_{N,r} \le L^2$. 


\subsubsection{Eigenvalue computations with finite difference method (bounded domain)}

The operator $L_a(\lambda)$ (see \Ref{FJeq16}) is defined on the domain 
$D({\cal A}) = \left\{ u \in H^2(\BigR) , x^4 u \in L^2 (\BigR) \right\}$. 
So $u$ is decreasing when $x^4$ is increasing and the decrease is faster than for the rotated harmonic oscillator \Ref{FJeq2}. So we want to consider the following nonlinear eigenvalue problem in bounded domain with Dirichlet homogeneous boundary conditions~: find $\lambda \in \BigC$ such that~:
\begin{equation}
\left\{
\begin{array}{l}
L_a(\lambda) u = 0 \ \ , \ \ x \in \Omega \\
u(\pm L) = 0 
\end{array}
\right.
\label{FJeq21}
\end{equation}
where $\Omega = (-L,+L)$ with $L$ sufficiently large. More precisely we retain $L = \sqrt{2N-2}$ (see \Ref{eq6}).  \\

As before, the problem \Ref{FJeq21} car be reformulated as an eigenvalue problem~: 
\begin{equation}
{\mathcal  A}_a \left(
\begin{array}{l} 
u \\
v
\end{array}
\right) 
= \lambda 
\left(
\begin{array}{l}
u \\
v
\end{array}
\right) 
\label{FJeq22}
\end{equation}
where~:
$$
                      {\mathcal  A}_{a} = \left(
                            \begin{array}{cc}
                            0 & \1\\
                            -L_0&-L_{1_a}\\
                            \end{array}
                            \right)
$$
with $v = \lambda u$ and the operators $L_0 = - \Frac{d^2}{dx^2} + x^4$, $L_{1_a} = - 2 a x^2$. \\

We consider on the domain $\Omega$ a meshgrid with a mesh $\Delta x = 2L/N$ on $\Omega$ and we note $x_j = -L + j \Delta x$, $j=0,\ldots N$ the points of the grid. We have retained homogeneous Dirichlet boundary conditions for $x=\pm L$, so $u(x_0) = u(x_N) =0$. We look for an approximation $u_N$, $v_N$ of $u$ and $v=\lambda u$ such that~:
$$
{\mathcal A}_{a,N}
 \left(
\begin{array}{l} 
U_N \\
V_N
\end{array}
\right) 
= \lambda_N 
\left(
\begin{array}{l}
U_N \\
V_N
\end{array}
\right) 
$$
with $U_N$ and $V_N$ two vectors containing respectively the approximations $u_N(x_j)$, $v_N(x_j)$ of $u(x_j)$, $v(x_j)$ and ${\mathcal A}_{a,N}$ is the square matrix of order $2N-2$~:
\beq\label{eq110}
                      {\mathcal  A}_{a,N} = \left(
                            \begin{array}{cc}
                            0 & \1_N\\
                            -L_{0_N}&-L_{1_{a,N}}\\
                            \end{array}
                            \right)
\eeq
where $L_{0_N} u_N (x_j) = -(\Frac{u_N(x_{j+1})-2u_N(x_j)+u_N(x_{j-1})}{\Delta x^2} + x_j^4 u_N(x_j)$ is the discretization of the operator $L_0$ with a centered finite difference scheme and $L_{1_{a,N}} v_N (x_j) = -2 a x_j^2 v_N(x_j)$. \\

$L_{0_N}$ is a tridiagonal symmetric matrix such that $L_{0_N}(j,j) = \Frac{2}{\Delta x^2} + x_j^4$ and $L_{0_N}(j,j-1) = - \Frac{1}{\Delta x^2}$.\\

$L_{1_{a,N}}$ is a diagonal matrix such that $L_{1_{a,N}}(j,j) = -2 a x_j^2$.\\

For the numerical computation of the spectrum of the matrix ${\mathcal A}_{a,N}$ we use the function DGEEV of the LAPack library.\\

Now, we are interested to analyze the dependence of the spectrum of the operator \Ref{FJeq16} in function of the real parameter $a$. For this, we consider an approximation of the infinite dimensional domain as a bounded domain with periodic boundary conditions. We look for eigenfunction $u_k (x) = \hat{u}_k \exp (i k' x)$, with $k'=\Frac{k \pi }{L}$, of the continuous operator \Ref{FJeq16}. Computing $L_a(\lambda) u_k(x)$ we obtain the following equation~:
$$
\lambda^2 - 2 a \lambda x^2 + x^4 + k'^2 = 0
$$
The discriminant $\Delta = 4 (a^2 - 1 ) x^4 - 4 k'^2$ is negative for $0 \le a \le 1$ and the solutions are~:
$$
\lambda^\pm = \Frac{2 a x^2 \pm i \sqrt{-\Delta}}{2}
$$
When $a$ is increased from $0$ to $1$ the ratio of the imaginary part over the real part of $\lambda$, $\Frac{\vert \lambda_i\vert }{\vert \lambda_r\vert } = \Frac{\sqrt{k'^2-(a^2-1) x^4 }}{a x^2}$ is decreased and it is infinite for $a=0$. We can observed this on the numerical simulations corresponding to $a=0$, $a=0.5$, $a=0.9$ and $a=1.0$ obtained with the finite difference scheme for $N=50$ and $L=10$ (see Figure \ref{Figure9} ). \\

We are now interested with the operator $L_a(\lambda)$ for $a=1$. 
We look for the spectrum of the discretized operator, using finite difference method, where $x$ is replaced with a real constant $b$ (see \Ref{FJeq17ter}). We have~:
\begin{equation}
\lambda_N^2 u(x_j) - \Frac{u(x_{j+1}) - 2 u(x_{j}) + u(x_{j-1})}{\Delta x^2} - 2 b^2 \lambda_N u(x_j) + b^4 u(x_j) = 0
\label{eq100}
\end{equation}
If we consider periodic boundary conditions, we look for a solution of \Ref{eq100} of the form $u(x) = \hat{u}_k \exp(ik'x)$, with $k'=\Frac{k\pi}{L}$. Substituting in \Ref{eq100} and supposing that $\hat{u}_k \ne 0$ we obtain~:
$$
\lambda_N^2 \Delta x^2 - 2 b^2 \lambda_N \Delta x^2 + b^4 \Delta x^2 - 2 \cos(k'\Delta x) + 2 = 0
$$
Finally we have $\lambda_N = \lambda_{N,r} + i \lambda_{N,i}$ with $\lambda_{N,r} = b^2$ is wavenumber independent and $\lambda_{N,i} = \pm \Frac{\sqrt{2 - 2 \cos (k'\delta x)}}{\Delta x}$ is wavenumber dependent. So the spectrum of the discretized operator is located in the part of the plan complex such that $\lambda_{N,r} = b^2$ and $-\vert k'_{\rm max} \vert \le \lambda_{N,i} \le \vert k'_{\rm max} \vert$ since $\cos(k'\Delta x) \simeq 1-\Frac{k'^2 \Delta x^2}{2}$ for $\Delta x$ sufficiently small.  \\

Here since $\Omega = (-L,+L)$ and $N$ is the number of grid points retained, the highest wavenumber $k'_{\rm max}$ we can take into account with this meshgrid is $k'_{\rm max} = \Frac{N}{2L}=\Delta x^{-1}$. Since $L \simeq \sqrt{2N-2}$ (see \Ref{eq6}) we have $k'_{\rm max} = O(\sqrt{N})$, which is in agreement with the Hermite spectral method for unbounded domain (see \Ref{FJeq20}). \\

On Figure \ref{Figure10} we present the spectrum of the matrix \Ref{eq110} for $N=50$, $L=10$ and $a=1$. Comparison with Figure \ref{Figure8}  shows that the results obtained for Hermite spectral method (unbounded domain) and for finite difference method (bounded domain) are quite similar. We have chosen $L=10$ for the size of the bounded domain, in agreement with \Ref{eq6}. As it has been said previously, theoretical results give that the eigenvalues of the continuous operator \Ref{FJeq16}, for $a=1$, are included in the two  sectors $\{\lambda\in\C,\; \vert\arg(\lambda)\vert \geq \Frac{\pi}{3}\}$ (see Section 3). But we can see on Figure \ref{Figure10} (as on Figure \ref{Figure8}) that computed eigenvalues are not all included in these two sectors, which can be imputed to numerical instabilities leading to spurious eigenvalues (spectral pollution, see \cite{DaviesPlum}). \\

In order to analyze these numerical instabilities, we study the stability of the eigenvalues in function of a perturbation on the points of the mesh grid retained for the discretization. 
The equality \Ref{eq14} measures the sensivity of the eigenvalue $\lambda_N$ of the matrix ${\mathcal A}_N$ in function of a perturbation $\varepsilon$ on the meshgrid (condition number of the eigenvalue $\lambda_N$). Here the matrix ${\mathcal E}_N$ is the matrix of order $2N-2$~:
$$
                      {\mathcal E}_N = \varepsilon \left(
                            \begin{array}{cc}
                            0 & 0 \\
                            {\mathcal E}_{0,N}&{\mathcal E}_{1,N}\\
                            \end{array}
                            \right).
$$
where ${\mathcal E}_{0,N}$ (resp. ${\mathcal E}_{1,N}$) is the diagonal matrix with the elements $-4 x_j^3$ (resp. $4 a x_j$) on the diagonal, $j=1,\ldots,N-1$ (we have neglected in ${\mathcal E}_N$ the terms in $\varepsilon^n$, with $n >1$). \\

On Figure \ref{Figure11} we have represented the condition number of the eigenvalues $\lambda_N$ in function of the modulus $\vert \lambda_N \vert $ for $N=50$, $L=10$ and $a=1$. We can see that eigenvalues are ill conditioned, excepted for the eigenvalues with small modulus. This can explain the convergence problem when $N$ is increased. In comparison with the rotated harmonic oscillator (see Figure \ref{Figure6}) we can see that the condition numbers of the eigenvalues are much greater for the operator \Ref{FJeq16} than for the rotated harmonic oscillator \Ref{FJeq2}. A small perturbation on the grid points induces large perturbations on the eigenvalue computations. However the eigenvalues are independent of $x$. So, in order to decrease this dependence of the eigenvalues in function of a perturbation on the points of the mesh grid, we have considered several grids for the finite difference discretization, with a shift on the mesh points, but with the same step $\Delta x$ for the mesh grid~: $y_j = x_j + \varepsilon$. Then we compute an average on the eigenvalues obtained with these staggered grids. The results obtained are presented on Figure \ref{Figure12}, which corresponds to $a=1$, $N=1000$ and $L=10$. The number of staggered grids retained is $11$. We can see that spurious eigenvalues have disappeared. The computed eigenvalues after averaging are now essentially contained in the area $\{\lambda\in\C,\; \vert\arg(\lambda)\vert \geq \Frac{\pi}{3}\}$ in agreement with theoretical results (see Section 3). We can note on Figure \ref{Figure12} that on the imaginary axis we have limited the imaginary part of $\lambda_N$ to $\vert \lambda_{N,i} \vert \le \Frac{N}{2L} = 50$. Indeed, as its has been said previously, $\lambda_{N,i}$ is function of the wavenumber and the highest wavenumber we can take into account on the grids is $\Frac{N}{2L}$. 

\noindent For the use of staggered meshes to avoid spectral pollution, we may mentioned the following reference \cite{LABV}. \\

Now we consider the pseudospectra \Ref{eq15} since it is known that the numerical computation of the pseudospectra is more stable than for the spectra (see Section 4). For the computation of the pseudospectra, we have retained complex values $z$ in \Ref{eq15} lying on the meshgrid in the part of the complex plane corresponding to $[0,100] \times [0,100]$. The step retained is $dx=dy=1$ in the real and imaginary directions. On Figure \ref{Figure13} we can see the computation for the matrix ${\mathcal A}_{a,N} - z \1_N$ with ${\mathcal A}_{a,N}$ corresponding to the matrix \Ref{eq110}, for $N=1000$, $L=10$ and $a=1$. We can note that, in agreement with the theoretical results (see Section 3), the two sectors $\{\lambda\in\C,\; \vert\arg(\lambda)\vert \geq \Frac{\pi}{3}\}$ of the spectrum of the continuous operator \Ref{FJeq16} are essentially contained in the area of the pseudospectra corresponding to the smallest values of the parameter $\varepsilon$, {\it i.e.} in the area where the distance of $z$ to the eigenvalues of the matrix \Ref{eq110} is the smallest. \\

The pseudospectra computation is very expensive. So we use parallel computation in order to accelerate the computation. 
The numerical solution is done thanks to the linear algebra library LAPack which contains specialized algorithms for singular values problems, especially the one called ZGESVD for complex matrices in double precision. As the matrix \Ref{eq110} is quite huge, and computing time a bit long, a parallelization by MPI (Message Passing Interface) is implemented with the client/server model. One process (the server) distributes values of the complex parameter $z$ (see \Ref{eq15}) to the other processes (the clients) which sample the domain. The server renews their data as the work progresses. Each client builds the matrix to be study and sends to the server, at the end of the computation, the smallest value. This system has the advantage of being dynamically balanced. As there is no communication (in MPI sense) between the clients, the efficiency of the parallelization is complete.
As an example, the simulation corresponding to the parameters $N= 5000$, $L=1000$, $a=1$ and to an area of the complex plane $[0,150]\times [0,150]$ with a mesh step $dx=1$ and $dy=1$ in the real and imaginary directions has needed 40 cores (Intel Xeon E5-2670 at 2.5GHz) during quite 40 days. \\

Now, as for the rotated harmonic oscillator, we consider here the computation of the pseudospectra based on Definition \ref{Def2} (see \Ref{Pseudo1}) instead of Definition \ref{Def1} as previously. So we look for $z \in \BigC$ such that~:
\begin{equation}
\vert \vert {A}^{-1}_{a,N}(z)  \vert\vert =  s^{-1}_{\rm min} ({A}_{a,N} (z) )\ge \varepsilon^{-1} \exp(\delta \vert z \vert ^\mu)
\label{Pseudo3}
\end{equation}
where ${A}_{a,N} (z)$ is the matrix obtained with the finite difference discretization of the operator $L_{a} (z)$ (see \Ref{FJeq16}); $\varepsilon$, $\delta$ and $\mu$ are real parameters. The majoration in \Ref{Pseudo3} depends on $\vert z \vert$, in opposition with \Ref{eq15}, {\it i.e.} the constraint appearing in \Ref{Pseudo3} is stronger when $\vert z \vert$ is increased. 
As before we have retained $N=1000$, $L=10$, $a=1$ and complex values $z$ lying on the meshgrid in the area of the complex plane corresponding to $[0,100] \times [0,100]$. The step retained is $dx=1$ and $dy=1$ in the real and imaginary directions. 
In order to look for the influence of the parameters $\delta$ and $\mu$ on the pseudospectra \Ref{Pseudo3}, we have presented on Figure \ref{Figure101} the pseudospectra computed with different values of the parameters $\delta$ and $\mu$. We can see that when the parameter $\mu$ is increased, eigenvalues with large modulus are eliminated in the pseudospectra computed with \Ref{Pseudo3}. Moreover, the CPU time required to compute pseudospectra with Definition \ref{Def2} (see \Ref{Pseudo3}) is much lower than if we use Definition \ref{Def1} (see \Ref{eq15}). Indeed, the matrix ${A}_{a,N} (z)$ is of order $N+1$ instead of $2N+2$ for the matrix ${\mathcal A}_{a,N}$. 

\subsubsection{Eigenvalue computations with Legendre spectral Galerkin method (bounded domain)}

In order to obtain a higher accurate numerical scheme in bounded domain, we propose a spectral numerical scheme using Legendre Galerkin basis.  \\

We consider the problem \Ref{FJeq21}. This problem is reformulated as an eigenvalue problem \Ref{FJeq22}. But instead of using a finite difference scheme to obtain an approximation $u_N$, $v_N$ of $u$ and $v=\lambda u$, we use a spectral method with Legendre Galerkin basis $\Phi_l$. Such basis is obtained as a linear combination of Legendre polynomials~:
$$
\Phi_l (x) = c_l (L_l(x)-L_{l+2}(x))
$$
with $L_l$ the Legendre polynomial of degree $l$ and $c_l = \Frac{1}{\sqrt{4l+6}}$ (see \cite{Jie}). Such a basis verify homogeneous Dirichlet boundary conditions $\Phi_l(\pm 1) = 0$. In particular, with the scalar product in $L_2(\Omega)$ we have~:
\begin{equation}
(\Phi_k,\Phi_j) = \left\{
\begin{array}{l}
c_k c_j ( \Frac{2}{2j+1} + \Frac{2}{2j+5} ) \ \ , \ \ k=j \\
-c_k c_j \Frac{2}{2k+1} \ \ , \ \ k=j+2 \\
0 \ \ , \ \ otherwise
\end{array}
\right.
\label{eq200}
\end{equation}
and
\begin{equation}
(\Phi'_k,\Phi'_j) = \left\{
\begin{array}{l}
1 \ \ , \ \ k=j \\
0 \ \ , \ \ k \ne j
\end{array}
\right.
\label{eq201}
\end{equation}
Moreover, we need the expressions of $x^2 \Phi_l$ and $x^4 \Phi_l$ as linear combination of the Legendre polynomials. We have~:
\begin{equation}
x^2 L_l(x) = \Frac{1}{2l+1} \left( \Frac{l+1}{2l+3} ((l+2) L_{l+2}(x) + (l+1) L_l(x)) + \Frac{l}{2l-1}  (l L_l(x) + (l-1)L_{l-2}(x))\right)
\label{eq201ter}
\end{equation}
and
\begin{equation}
x^4 L_l(x) = \alpha_l L_{l-4}(x) + \beta_l L_{l-2} (x) + \gamma_l L_l(x) + \delta_l L_{l+2} (x) + \eta_l L_{l+4}(x) 
\label{eq201qua}
\end{equation}
with
$$
\alpha_l = \Frac{1}{2l+1} \left( \Frac{l(l-1)(l-2)(l-3)}{(2l-1)(2l-3)(2l-5)} \right)
$$
$$
\beta_l = \Frac{1}{2l+1} \left( (\Frac{(l+1)^2}{(2l+3)(2l+1)} + \Frac{l^2}{(2l-1)(2l+1)})(\Frac{l(l-1)}{2l-1}) + \Frac{l(l-1)^3}{(2l-1)^2(2l-3)} + \Frac{l(l-1)(l-2)^2}{(2l-1)(2l-3)(2l-5)} \right)
$$
$$
\gamma_l = \Frac{1}{2l+1} \left( (\Frac{(l+1)^2(l+2)^2}{(2l+3)^2(2l+5)} + (\Frac{(l+1)^2}{(2l+1)(2l+3)} + \Frac{l^2}{(2l-1)(2l+1)}) * (\Frac{(l+1)^2}{(2l+3)} + \Frac{l^2}{(2l-1)}) + \Frac{l^2(l-1)^2}{(2l-1)^2(2l-3)} \right)
$$
$$
\delta_l = \Frac{1}{2l+1} \left( (\Frac{(l+1)^2}{(2l+3)(2l+1)} + \Frac{l^2}{(2l-1)(2l+1)})(\Frac{(l+1)(l+2)}{2l+3}) + \Frac{(l+1)(l+2)^3}{(2l+3)^2(2l+5)} + \Frac{(l+1)(l+2)(l+3)^2}{(2l+3)(2l+5)(2l+7)} \right)
$$
$$
\eta_l = \Frac{1}{2l+1} \left(\Frac{(l+1)(l+2)(l+3)(l+4)}{(2l+3)(2l+5)(2l+7)} \right)
$$
In order to adapt the previous basis $\Phi_l$ to the Dirichlet boundary conditions $\Phi_l(\pm L) = 0$, we multiply the previous polynomials by a scale factor. As for the Hermite spectral method (see Section 5.2.1), we use a method of weighted residuals (MWR, see for example, \cite{GO}, \cite{Cal}) and relations \Ref{eq200}, \Ref{eq201}, \Ref{eq201ter}, \Ref{eq201qua} to obtain the following generalized eigenvalue problem~:
\begin{equation}
{\mathcal A}_{a,N}
 \left(
\begin{array}{l} 
U_N \\
V_N
\end{array}
\right) 
= \lambda_N {\mathcal B}_N 
\left(
\begin{array}{l}
U_N \\
V_N
\end{array}
\right) 
\label{eq201bis}
\end{equation}
where $U_N$ and $V_N$ are the vectors containing respectively the coefficients $\tilde{u}_l$ and $\tilde{v}_l$, $l=0,\ldots,N$, of $u_N=\sum_{l=0}^N {\tilde u}_l \Phi_l$ and $v_N=\sum_{l=0}^N {\tilde v}_l \Phi_l$. ${\mathcal A}_{a,N}$ is the square matrix of order $(2N+2)$~:
$$
                      {\mathcal A}_{a,N} = \left(
                            \begin{array}{cc}
                            0 & \1_N\\
                            -L'_{0_N}&-L'_{1_{a,N}}\\
                            \end{array}
                            \right)
$$
and ${\mathcal B}_N$ is the square matrix of order $(2N+2)$~:
$$
                      {\mathcal B}_N = \left(
                            \begin{array}{cc}
                            B_{0_N} & 0\\
                            0 & B_{0_N}\\
                            \end{array}
                            \right)
$$
Here $L'_{0_N} u_N = (L_0 u_N,\Phi_{l'})$ and $L'_{1_{a,N}} v_N = (L_{1_a} v_N,\Phi_{l'})$, $l'=0,\ldots,N$, with
$L_0 u_N = - \Frac{d^2u_N}{dx^2} + x^4 u_N$ and $L_{1_a} u_N = - 2 a  x^2 u_N$. As for $B_{0_N} u_N = (u_N,\Phi_{l'})$. \\

$L'_{0_N}$ is a symmetric matrix with seven diagonal and $L'_{1_{a,N}}$ is a pentadiagonal symmetric matrix. As for $B_{0_N}=(\Phi_l,\Phi_{l'})$ for $l$ and $l'=0,\ldots N$ (see \Ref{eq200}). \\ 

To obtain the eigenvalues of the generalized eigenvalue problem \Ref{eq201bis} we use the function DGGEV of the LAPack library.\\

On Figure \ref{Figure14} we present the solutions $\lambda_N$ of \Ref{eq201bis}, computed with $N=50$, $L=10$ and $a=1$. Comparison with the spectral Hermite method (Figure \ref{Figure8}) and the finite difference method (Figure \ref{Figure10}) is done. We can see that the numerical results are quite similar. 

\subsection{Another discretized nonlinear eigenvalue problem}

In this section we consider the following operator~:
\begin{equation}
\mathcal{L} u(x) =-\Frac{d^2u}{dx^2} (x)+ (x^{k}-\lambda)^2u(x) \ 
\label{adnlep}
\end{equation}
For $k=2$ we retrieve the operator \Ref{Deq3} studied in the previous section.  \\

We discretize the problem $\mathcal{L} u=0$ using some techniques similar to finite difference methods, with a spatial step equal to one. For simplicity reasons we need to add either periodic boundary conditions or homogeneous boundary conditions. Also we replace $\Delta u(n)$ by $\delta \delta^{*}$ where:
$$
\begin{array}{c}
\delta u(n) = u(n+1) - u(n),  \ n \in \N \\
\delta^{*} u(n)= u(n) - u(n-1),  \ n\in \N
\end{array}
$$
 i.e.
$$
(\delta \delta^{*}) u(n) = u(n-1)  - 2 u(n) + u(n+1)
$$
 So we have :
\begin{equation}\label{disprobFD}
\mathcal{L} u(n) =- (\delta \delta^{*})u(n)+ (n^{k}-\lambda)^2 u(n),  \quad n \in \N
\end{equation}
\subsubsection{Finite difference method with periodic boundary conditions}

In this section we are interested to study the problem (\ref{disprobFD}) with periodic boundary conditions. So for some $N\in\N$, we study the following problem :
$$
\begin{array}{c}
-(\delta \delta^{*}) u(n) + (n^{k} -\lambda )^2  u(n) =  0, \quad n=1,\cdots,N \\
u(j) = u(j+N), \quad j = 0,1
\end{array}
$$
For $n=1,\cdots, N$ we have :
$$
\begin{array}{ccccc}
\mbox{\small{ n=1}}   \hfill:& \hfill -u(0)+2u(1)-u(2)+1^{2k} u(1)-2\lambda (1)^k u(1) +\lambda^2 u(1)&=&0\\
\mbox{\small{ n=2}}  \hfill:& \hfill -u(1)+2u(2)-u(3)+2^{2k} u(2)-2\lambda (2)^k u(2) +\lambda^2 u(2)&=&0\\
\mbox{\small{ n=j}} \hfill:& \hfill -u(j-1)+2u(j)-u(j+1)+j^{2k} u(j)-2\lambda (j)^k u(j) +\lambda^2 u(j)&=&0\\
\mbox{\small{ n=N-1}} \hfill:&  \hfill  -u(N-2)+2u(N-1)-u(N)+(N-1)^{2k} u(N-1) \\
             &\hfill -2\lambda (N-1)^k u(N-1)+\lambda^2 u(N-1)&=&0\\
\mbox{\small{ n=N}}   \hfill:& \hfill -u(N-1)+2u(N)-u(N+1)+(N)^{2k} u(N)-2\lambda (N)^k u(N) +\lambda^2 u(N)&=&0\\
\end{array} 
$$
Using the periodic conditions $u(0)=u(N)$ and $u(N+1)=u(1)$, we obtain the system~:
$$
\begin{array}{ccccc}
\mbox{\small{ n=1}}   \hfill:& \hfill -u(N)+2u(1)-u(2)+1^{2k} u(1)-2\lambda (1)^k u(1) +\lambda^2 u(1)&=&0\\
\mbox{\small{ n=2}}  \hfill:& \hfill -u(1)+2u(2)-u(3)+2^{2k} u(2)-2\lambda (2)^k u(2) +\lambda^2 u(2)&=&0\\
\mbox{\small{ n=j}} \hfill:& \hfill -u(j-1)+2u(j)-u(j+1)+j^{2k} u(j)-2\lambda (j)^k u(j) +\lambda^2 u(j)&=&0\\
\mbox{\small{ n=N-1}} \hfill:&  \hfill  -u(N-2)+2u(N-1)-u(N)+(N-1)^{2k} u(N-1) \\
             &\hfill -2\lambda (N-1)^k u(N-1)+\lambda^2 u(N-1)&=&0\\
\mbox{\small{ n=N}}   \hfill:& \hfill -u(N-1)+2u(N)-u(1)+(N)^{2k} u(N)-2\lambda (N)^k u(N) +\lambda^2 u(N)&=&0\\
\end{array} 
$$
This gives the following system :
$$A_0+\lambda A_1+\lambda^2 I=0$$
where $I$ is the $N\times N$ identity matrix and $A_1$, $A_0$ are given as follows :
\begin{equation} \label{A1}
A_1 = -2 \left( \begin{array}{ccccccc}

          1 & 0     &\cdots &  &  & \cdots & 0 \\

          0 & 2^k   &0   & \cdots &  & \cdots & 0 \\

 \vdots  &   & \vdots  &  & \vdots &  &   \vdots  \\
 0 & \cdots &  &\cdots  &0  &       (N-1)^k &0 \\
          0 & \cdots &  & &\cdots    &        0 &N^k \\
\end{array} \right)
\end{equation} 
\begin{equation}\label{A0} 
A_0=A_{0,d}+A_{0,+1}+A_{0,-1}
\end{equation} 
with 
$$
A_{0,d} =\left( \begin{array}{cccccccccc}
          2+1 & 0    & 0  & \cdots &  &         &  & \cdots &0  & 0\\     
           0 & 2+2^{2k}  & 0 & 0 &\cdots       &  &  & & \cdots   & 0  \\
              & &        &\vdots &  &  & \vdots &  &  & \\
           0  & \cdots &  \cdots     &0 &  0    & 2+j^{2k} & 0 &  0  & \cdots & 0  \\
             \vdots   &  &  &     \vdots    &  &  & \vdots &  &  &  \vdots \\
           0  & \cdots &        &      &  & \cdots & 0 & 0& 2 + (N-1)^{2k} &  0   \\
          0  &0& \cdots&  &   &        &    \cdots   & 0 & 0         & 2 + N^{2k}   \\
\end{array} \right) 
$$
$$
A_{0,+1} = \left( \begin{array}{ccccccc}

          0 & -1     &\cdots &  &  & \cdots & -1 \\

          0 & 0   &-1  & \cdots &  & \cdots & 0 \\

 \vdots  &   & \vdots  &  & \vdots &  &   \vdots  \\
 0 & \cdots &  &\cdots  &0  &     0 &-1 \\
          0 & \cdots &  & &\cdots    &        0 &0 \\
\end{array} \right), \quad 
A_{0,-1} = \left( \begin{array}{ccccccc}

          0 & 0    &\cdots &  &  & \cdots & 0 \\

          -1 & 0   &0 & \cdots &  & \cdots & 0 \\
0 & -1  &0 &\cdots   &    & \cdots & 0 \\
 \vdots  &   & \vdots  &  & \vdots &  &   \vdots  \\
0 & \cdots &  &\cdots  &0  &     0 &0 \\
          -1 & \cdots &  & &\cdots    &        -1 &0 \\
\end{array} \right)
$$
%
%
We start computing the eigenvalues for different values of $N$ and for the operator  $\mathcal{L}$. Then, we compute the eigenvalues for some perturbations of the operator $\mathcal{L}$, i.e. we study the discrete operator :
$$\mathcal{L}_{c}u(n)=A_0 u(n)+c\lambda A_1 u(n) + \lambda^2 I u(n) , \quad n=1,\cdots, N $$
for $ 0\leq c\leq 1$ with the same previous periodic boundary conditions. For this we consider the linearization system problem in place of the non-linear problem, so we study the spectrum of the linear system $\mathcal{A}_c U=\lambda U$ with :
$$\mathcal{A}_c=\left( \begin{array}{cc} 
               0 & I \\
               -A_0&-cA_1
\end{array} \right),  \quad 0\leq c\leq 1 $$
where $U=(u_1,u_2,\cdots,u_{N-1},u_N,v_1,v_2,\cdots,v_{N-1},v_N)^t$, with $v_i = \lambda u_i$, $i=1,\cdots,N$. $A_0$ and $A_1$ are given in (\ref{A0}) and (\ref{A1}) respectively.
For the computation of the eigenvalues, we use Matlab (or Scilab). \\

The results obtained for $N=100$ $k=2$ and $c=1$ are presented on Figure \ref{zomec1}. The associated domain is $[0,N]$. This figure represents a zoom for the case $c=1$. We note that the imaginary part of the eigenvalues $\lambda_i $ lies between $1.38$ and $1.42$ in the positive part and between $-1.42$ and $-1.38$ in the negative part. Starting from a real part $\lambda_r= 576$ all the eigenvalues are aligned on a straight parallel to the $x-axis$ with $\lambda_i=1.4141$ and $\lambda_i=-1.4141$. The results obtained for $N=1000$ $k=2$ and $c=1$ are similar. \\

On Figure \ref{les6cassimplex4} we present the numerical results obtained for $N=1000$, $k=4$ and $0\leq c \leq1$. For the case $c=0$ we have pure imaginary eigenvalues (since in this case we have just a selfadjoint matrix). The positions of eigenvalues for the cases $c= 0.2$, $0.4$ confirm the theoretical results. For the cases $c=0.6$, $0.8$, $1$, eigenvalues are localized in a sector delimited by an angle with the $x-axis$ smaller than $2\pi/6$. This is not coherent with the theoretical results.

\subsubsection{Finite difference method with homogeneous boundary conditions}
   In the following we consider the problem (\ref{disprobFD}) with homogeneous boundary conditions. So we study the following problem :
\begin{equation}
\begin{array}{c}
-(\delta \delta^{*}) u(n) + (n^{k} -\lambda )^2  u(n)  =  0, \quad n=1,\ldots,N \\
u(0)=u(N+1)=0
%
%
\end{array}
\label{cask6}
\end{equation}
So we obtain the following system~:
$$A_0+\lambda A_1+\lambda^2 I=0$$
where $I$ is the $N\times N$ identity matrix and $A_1$, $A_0$ are given as follows :
\begin{equation} \label{A1h}
A_1 = -2 \left( \begin{array}{ccccccc}

          1 & 0     &\cdots &  &  & \cdots & 0 \\

          0 & 2^k   &0   & \cdots &  & \cdots & 0 \\

 \vdots  &   & \vdots  &  & \vdots &  &   \vdots  \\
 0 & \cdots &  &\cdots  &0  &       (N-1)^k &0 \\
          0 & \cdots &  & &\cdots    &        0 &N^k \\
\end{array} \right)
\end{equation} 
\begin{equation}\label{A0h}
A_0=A_{0,d}+A_{0,+1}+A_{0,-1}
\end{equation} 
where
$$
A_{0,d} =\left( \begin{array}{cccccccccc}
          2+1 & 0    & 0  & \cdots &  &         &  & \cdots &0  & 0\\     
           0 & 2+2^{2k}  & 0 & 0 &\cdots       &  &  & & \cdots   & 0  \\
              & &        &\vdots &  &  & \vdots &  &  & \\
           0  & \cdots &  \cdots     &0 &  0    & 2+j^{2k} & 0 &  0  & \cdots & 0  \\
             \vdots   &  &  &     \vdots    &  &  & \vdots &  &  &  \vdots \\
           0  & \cdots &        &      &  & \cdots & 0 & 0& 2 + (N-1)^{2k} &  0   \\
          0  &0& \cdots&  &   &        &    \cdots   & 0 & 0         & 2 + N^{2k}   \\
\end{array} \right) 
$$
$$
A_{0,+1} = \left( \begin{array}{ccccccc}

          0 & -1     &\cdots &  &  & \cdots & 0 \\

          0 & 0   &-1  & \cdots &  & \cdots & 0 \\

 \vdots  &   & \vdots  &  & \vdots &  &   \vdots  \\
 0 & \cdots &  &\cdots  &0  &     0 &-1 \\
          0 & \cdots &  & &\cdots    &        0 &0 \\
\end{array} \right), \quad
A_{0,-1} = \left( \begin{array}{ccccccc}

          0 & 0    &\cdots &  &  & \cdots & 0 \\

          -1 & 0   &0 & \cdots &  & \cdots & 0 \\
0 & -1  &0 &\cdots   &    & \cdots & 0 \\
 \vdots  &   & \vdots  &  & \vdots &  &   \vdots  \\
0 & \cdots &  &\cdots  &0  &     0 &0 \\
          0 & \cdots &  & &\cdots    &        -1 &0 \\
\end{array} \right)
$$
We start by computing the eigenvalues for different values of $N$ and for the operator $\mathcal{L}$. Then we compute the eigenvalues for some perturbations of the operator $\mathcal{L}$, i.e. we consider the discrete operator :
$$\mathcal{L}_{c}u(n)=A_0 u(n)+c\lambda A_1 u(n) + \lambda^2 I u(n) , \quad n=1,\cdots,N$$
with $ 0\leq c\leq 1$ and the same previous homogeneous boundary conditions. We do this considering the linearization system problem in place of the non-linear problem. So we study the spectrum of the linear system $\mathcal{A}_c U=\lambda U$ with :
$$\mathcal{A}_c=\left( \begin{array}{cc} 
               0 & I \\
               -A_0&-cA_1
\end{array} \right),  \quad 0\leq c\leq 1 $$
where $U=(u_1,u_2,\cdots,u_{N-1},u_N,v_1,v_2,\cdots,v_{N-1},v_N)^t$, with $v_i = \lambda u_i$, $i=1,\cdots,N$. $A_0$ and $A_1$ are given in (\ref{A0h}) and (\ref{A1h}) respectively.
We compute the eigenvalues using Matlab. \\

For the numerical simulations we have considered a domain $[-L,+L]$ and a spatial step $\Delta x = \Frac{2L}{N}$. For the case $k=4$, the results obtained for $L=10$, $N=2000$ (resp. $L=20$, $N=10000$) and $c=1$ are presented on Figures \ref{N2000k4c1L10tt3fig} and \ref{N10000k4c1L20tt3fig} respectively. For the case $k=6$, the numerical results obtained for the example (\ref{cask6}) with $N=10000$, $c=1$ and $L=20$ (resp. $L=10$) are presented on Figure \ref{N10000k6c1L20tt3figs} and \ref{N10000k6c1L10tt3fig} respectively.

\begin{remark}
We can note that when the parameter $k$ is increased, the numerical results obtained are in better agreement with the theoretical results given in Section 3, {\it i.e} the eigenvalues of the continuous operator \Ref{adnlep} are included in the two  sectors $\{\lambda\in\C,\; \vert\arg(\lambda)\vert \geq \frac{k\pi}{2(k+1)}\}$. This can be explained by the fact that the eigenvalues are better conditionned when $k$ is increased.
\end{remark}

\section{Conclusions and open problems}

In this work we have presented a review of some theoretical results obtained for quadratic family of operators~:
$$
L(\lambda)=L_0+\lambda L_1 +\lambda^2
$$
where $L_0$ and $L_1$ are operators in an Hilbert space. \\

Then we have presented numerical methods to compute the spectrum of such operators. We reduce it to a  non self-adjoint linear eigenvalue problem. The numerical methods proposed are spectral methods and finite difference methods, for bounded and unbounded domains. For bounded domain we consider homogeneous Dirichlet boundary conditions and periodic boundary conditions. Comparison with the results obtained in unbounded and bounded domains are done. They are based on the size of the containment domain, deduces from the zeroes of the Hermite functions. \\

The numerical results obtained are presented. In particular the numerical instabilities are highlighted. 
Comparisons of the numerical results obtained, with the theoretical results presented in the first part of this work, are done. These comparisons show the difficulties for the numerical computation of such problem. Elimination of the spectral pollution, using staggered grids, and the computation of pseudospectra allow to obtain numerical results in agreement with theoretical results. \\

A future step in this work is the extension to the two dimensional case. This work is in progress and will be presented elsewhere. \\

\noindent{\bf Acknowledgments} \\

\noindent This work was initiated during the visit of Fatima Aboud, at the Laboratoire de Math\'ematiques Jean Leray,
Universit\'e de Nantes (France), CNRS UMR 6629. This visit was supported by the research project {\it D\'efiMaths} of
the F\'ed\'eration de Math\'ematiques des Pays de la Loire, CNRS FR 2962. Computations are done thanks to the computer of the CCIPL (Centre de Calcul Intensif des Pays de la Loire). \\

\begin{appendices}


\section{Hermite spectral  method}

\subsection{The  1-D  case}

 The basis $\{\varphi_k\}_{k\in\N}$ of Hermite functions is obtained as an orthonormal basis 
of $L^2(\R)$ of the eigenfunctions of the harmonic oscillator : 
   $$
   H_{osc} = -\Frac{d^2}{dx^2} + x^2
   $$
 We recall briefly its construction (see the basic books of quantum mechanics).\\
 Define the creation  operator $a^*$ and the annihilation operator $a$ 
$$
  a^* = x - \Frac{d}{dx} ,\;\; a = x + \Frac{d}{dx}
$$
We satisfy 
$$
 [a,a^*] = 2\1,\;\; H_{osc} = a^*a +\1 = \Frac{1}{2}(aa^* + a^*a)
$$
 where $ [a,a^*] =aa^* -a^*a$.
 
Starting by the normalized Gaussian :
$$
 \varphi_0(x) = \pi^{-1/4}{\rm e}^{-x^2/2} 
$$
verified $a\phi_0 = 0$ and then $H_{osc}\varphi_0 = \varphi_0$ one define by induction for integer $k$
 the sequence $\{\varphi_k\}_{k\in\N}$:
$$
 \varphi_{k+1} = (2(k+1))^{-1/2}a^*\varphi_{k}
$$
$$
\varphi_k = 2^{-k/2}(k!)^{-1/2}(a^*)^k\varphi_0
$$
We verify the following relation by using  an algebraic calculation
 \bea\label{crerec}
 a\varphi_{k+1} &=& (2k+1)^{1/2}\varphi_k \\
 a^*\varphi_{k} &=& (2k+1)^{1/2}\varphi_{k+1} \\
 H_{osc}\varphi_k &=& (2k+1)\varphi_k \\
 \la\varphi_k, \varphi_\ell\ra &=&\delta_{k,\ell}
 \eea
 where $\la, \ra$  denoted the scalar product in the (complex) Hilbert space $L^2(\R)$.\\
 We then show that $\{\varphi_k\}_{k\in\N}$ is a Hilbertian basis of $L^2(\R)$.
 
To do the projection of the  differential operators in this basis we need to calculate 
 the multiplication by $x$ and the derivation $\frac{d}{dx}$ of $\varphi_k$.\\
We use the relations $x = \frac{a + a^*}{2} $ and $\frac{d}{dx} = \frac{a - a^*}{2}$. By the
  relations (\ref{crerec}) we obtain : 
  \beq\label{psphix}
  \begin{array}{ccc}
   x\varphi_k &=& 2^{-1/2}\sqrt k\varphi_{k-1} + \sqrt{k+1}\varphi_{k+1} \hfill \\
   \frac{d}{dx}\varphi_k &=& 2^{-1/2}\sqrt k\varphi_{k-1} - \sqrt{k+1}\varphi_{k+1} \hfill\\
    x^2\varphi_k &=& \frac{1}{2}\left(\sqrt{k(k-1)}\varphi_{k-2} +(2k +1)\varphi_k +\sqrt{(k+1)(k+2)}\varphi_{k+2}\right)\hfill\\
     x^4 \varphi_k &= & \frac{1}{4}\left(\sqrt{k(k-1)(k-2)(k-3)}\varphi_{k-4} + (4k-2)\sqrt{k(k-1)}\varphi_{k-2} \hfill\right.\\
      & & +(6k^2+ 2k +3)\varphi_k + (4k +6)\sqrt{(k+1)(k+2)}\varphi_{k+2} \hfill \\   
      & &  + \left. \sqrt{(k+1)(k+2)(k+3)(k+4)}\varphi_{k+4}\right)\hfill\\
       \frac{d^2}{dx^2}\varphi_k &=& \frac{1}{2}\left(\sqrt{k(k-1)}\varphi_{k-2}  
        -  (2k+1)\varphi_k  + \sqrt{(k+1)(k+2)}\varphi_{k+2} \right) \hfill   
  \end{array} 
   \eeq
 We have used the following convention : when any integer become $<0$  we replace it by 0.
 
  \vspace{1cm}
 
 {\bf Estimation of error}
 
   \vspace{0.5cm}
   The suitable spaces are Sobolev spaces with weight are naturally associated to the harmonic oscillator
   $H_{osc}$ because the usual spaces of Sobolev are associated with the Laplacian. 
  For each integer $m\geq 0$ we define the space ${\cal B}_m$ of function $u\in L^2(\R)$ such that for any pair of integers $k,l$ such that $k+\ell \leq  m$ we have $x^k\frac{d^\ell}{dx^\ell}u\in L^2(\R)$.\\
   ${\cal B}_m$ is a Hilbert space with the scalar product
$$
 \la u, v\ra_m = \sum_{k+\ell \leq m}\int_\R \left(\overline{x^k\frac{d^\ell}{dx^\ell}}u\right)\left({x^k\frac{d^\ell}{dx^\ell}v}\right)dx
$$
     ${\cal B}_m$ is equal to the domain of $H_{osc}^{m/2}$ and the scalar product is equivalent to 
     $$
     \la u, v\ra^\star_m = \la H_{osc}^{m/2}u, H_{osc}^{m/2}v\ra = \la H_{osc}^mu, v\ra
     $$
 We deduce a characterization of  ${\cal B}_m$ with the  Hermite coefficient of $u$,
 $\alpha_k(u) :=  \la\varphi_k, u\ra$.
 \begin{proposition}\label{he}
 $u\in {\cal B}_m$ if and only if $\di{\sum_{k\in\N}(2k+1)^{m}\vert\alpha_k\vert^2 < +\infty}$.\\
  In addition, the scalar product is expressed as the following : 
 $$
   \la u, v\ra^\star_m = \sum_{k\in\N} (2k+1)^{m}\overline{\alpha_k(u)}{\alpha_k(v)}
   $$
   \end{proposition}
   The proposition can be summarized by saying that ${\cal B}_m$ is identical to the domain of the operator
   $H_{osc}^{m/2}$.  By complex interpolation we deduce the intermediate spaces ${\cal B}_s$
    for all $s$ positive reals hence by the duality for $s$ negative reals.  The arguments are identical to the case of usual Sobolev spaces. For $s<0$ the ${\cal B}_s$ are the spaces of temperate distribution. 
    
   Then we set $\di{u_N = \sum_{0\leq k\leq N}\alpha_k(u)\varphi_k}$, let $u_N=\Pi_Nu$,
   $\Pi_N$ be the projections on the vector space $V_N$  generated by $\{\varphi_0,\varphi_1,\cdots,\varphi_N\}$.
   So we clearly have : 
   $$
   \Vert u - u_N\Vert^{2}= \sum_{k >N}\vert\alpha_k\vert^2 \leq \frac{1}{(2N+1)^m}\sum_{k\in\N}(2k+1)^{m}\vert\alpha_k\vert^2
   $$
   Hence if  $u\in {\cal B}_m$ we have
   \beq\label{err1}
     \Vert u - u_N\Vert^2 \leq (2N+1)^{-m}\Vert u\Vert_m^{\star, 2}
  \eeq 
 More generally we can estimate the error in the spaces  ${\cal B}_s$
 \beq\label{err2}
     \Vert u - u_N\Vert_s^{\star,2} \leq (2N+1)^{s-m}\Vert u\Vert_m^{\star, 2}
  \eeq 
  It may be useful to have such Sobolev inequalities explaining the regularity and decay at infinity of
  $u \in {\cal B}_s$ as soon as $s$  is large enough.   We do not search to obtain an optimal estimation.
  \begin{proposition}
  Let $m\in \N$. There exists constants  $C_m>0$, $C_{s,m}$ ($m<2s -2$) such that
  \beq\label{hermest}
  \vert x^\ell\frac{d^k}{dx^k}\varphi_j(x)\vert \leq C_m (2j + 1)^{(m+1)/2},\;\; \forall x\in\R,\;\; k+\ell\leq m
  \eeq
  \beq\label{sobpo}
   \vert x^\ell\frac{d^k}{dx^k}u(x)\vert \leq C_{s,m} \Vert u\Vert_s,\;\; \forall x\in\R,\;\; k+\ell\leq m, u\in{\cal S}(\R)
   \eeq
   In particular if $m$ is known and if  $s> m +2$ then all $u$ in ${\cal B}_s$ are of class $C^m$ on $\R$
   and verify the inequality (\ref{sobpo}).
    \end{proposition}
{\em Proof}. 
It is sufficient to consider the case  $k=0$.\\
For (\ref{hermest})  starting from the usual Sobolev inequality (in one dimension the critical index 1/2).\\
Then
$$
\vert x^m\varphi_j(x)\vert \leq C \Vert x^\ell\varphi_j\Vert_{H^1} \leq C  \Vert\varphi_j\Vert_{m+1}
$$
Since $\Vert\varphi_j\Vert_m$ is of order $(2j+1)^{m/2}$, hence  (\ref{hermest}).\\
For (\ref{sobpo}) by expanding $u$ on the basis and by applying the Cauchy-Schwarz inequality
$$
\vert x^mu(x)\vert \leq \left(\sum_{j}\vert\alpha_j\vert^2(2j+1)^{s}\right)^{1/2}
 \left(\sum_{j}(2j+1)^{-s}\vert x^m\varphi_j(x)\vert^2\right)^{1/2}
 $$
then one use (\ref{hermest}) by choosing $s> m +2$
 $$
 \vert x^mu(x)\vert  \leq   C_{s,m} \Vert u\Vert_s
  $$
The last assertion follows from the density of ${\cal S}(\R)$ in ${\cal B}_s$ for all $s\in\R$.\\
The definition of derivatives does not cause a problem. \\
 Thus we see that the functions $u\in {\cal B}_m$  are both regular and decreasing to $0$
at the  infinity more rapidly when $m$ is big (positive).
 
 \hfill$\square$ 
\subsection{The  multidimensional  case}
The results are similar up to complication of notations.\\
Let $d\geq 2$, we denote $x=(x_1,x_2,\cdots, x_d)\in\R^d$, $k=(k_1,k_2,\cdots,k_d)\in\N^d$ and
$$
\varphi_k(x) = \varphi_{k_1}(x_1)\varphi_2(x_2)\cdots \varphi_d(x_d)
$$
Then $\{\varphi_k\}_{k\in\N^d}$ is an orthonormal basis of the Hilbert space 
 $L^2(\R^d)$.\\
 Then we have $d$  annihilation operators $a=(a_1,a_2,\cdots, a_d)$ and $d$ creation operators
   $a^*=(a^*_1,a^*_2,\cdots, a^*_d)$ with
   $$
     a_k^* = x_k - \frac{\partial}{\partial x_k} ,\;\; a_k = x + \frac{\partial}{\partial x_k}
     $$
To define the spaces ${\cal B}_m$ one replace $k, \ell$ by multi-indices.
The harmonic oscillator can be written as 
$$ H_{osc} = -\triangle  + \vert x\vert^2 = a^*\cdot a + d = \frac{1}{2}(a^*\cdot a + a\cdot a^*)
   $$
where $\di{a^*\cdot a = \sum_{1\leq k\leq d} a_k^* a_k}$ and $\vert x\vert^2 = x_1^2 +\cdots x_d^2$. \\
Then we have 
$$
H_{osc}\varphi_k = 2(k_1 + \cdots k_d) + 1,\; k = (k_1,\dots, k_d)
$$
The space $V_N$ is generated by $\{\varphi_k,\; k_1+\cdots k_d\leq N\}$. 
We denote for all multiindex $k$, $\vert k\vert = k_1 +\cdots k_d$.\\
The Hermite coefficients  
$\alpha_k(u)$ are indexed on $\N^d$. The estimation of error is then formally unchanged.\\
Sobolev inequalities with weight depend naturally on the dimension $d$.
For all $m, s$ such that $s > 1/2 + m+ d$ there exists $C_{s,m}>0$ such that
$$
 \vert x^ju(x)\vert  \leq   C_{s,m} \Vert u\Vert_s
  $$
For $\vert j\vert \leq m$ and  $u\in{\cal B}_s$.  Here $x^j = x_1^{j_1}\cdots x_d^{j_d}$ when
$j=(j_1,\cdots, j_d)$. \\
As for the case $d=1$ we have a similar inequality for the partial derivatives. 

\section{Figures}

In this section we give the figures referenced in this article.

\newpage
\phantom{aaaa}

\begin{figure}
  \centering
 \includegraphics[width=4in]{Figure100.eps}
\caption{ \label{Figure100} Computation of the pseudospectra \Ref{Pseudo2} of the matrix ${\mathcal A}_N$ corresponding to the finite difference scheme for $\alpha = \pi/4$, $N=100$ and $L=15$. }
\end{figure}
\begin{figure}
  \centering
\includegraphics[width=4in]{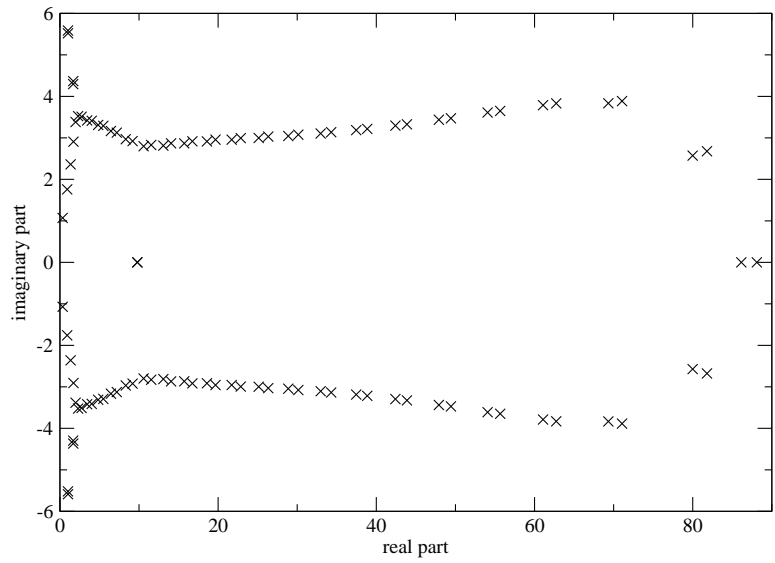}
\caption{ \label{Figure8} Spectrum of the matrix ${\mathcal A}_{a,N}$ \Ref{FJeq17bis} (Hermite spectral method) for $N=50$ and $a=1$. }
\end{figure}
\begin{figure}
  \centering
 \includegraphics[width=4in]{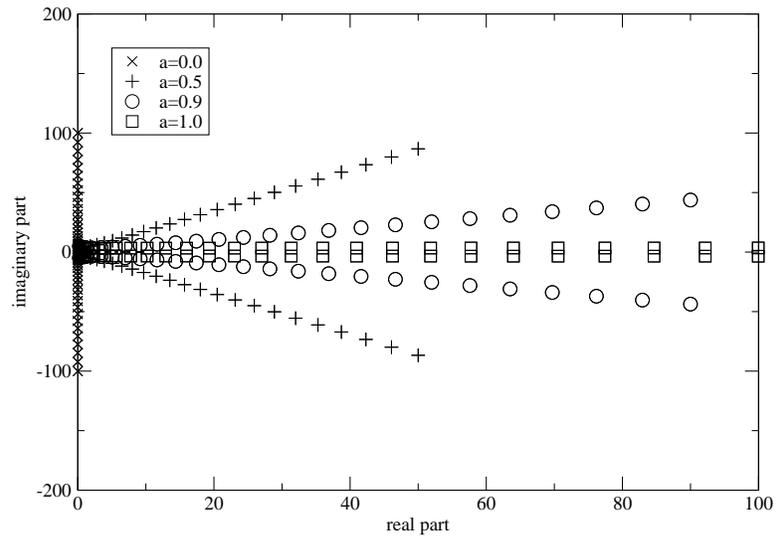}
\caption{ \label{Figure9} Spectrum obtained with the finite difference scheme for $a=0$, $a=0.5$, $a=0.9$, $a=1.0$, $N=50$ and $L=10$.}
\end{figure}
\begin{figure}
  \centering
 \includegraphics[width=4in]{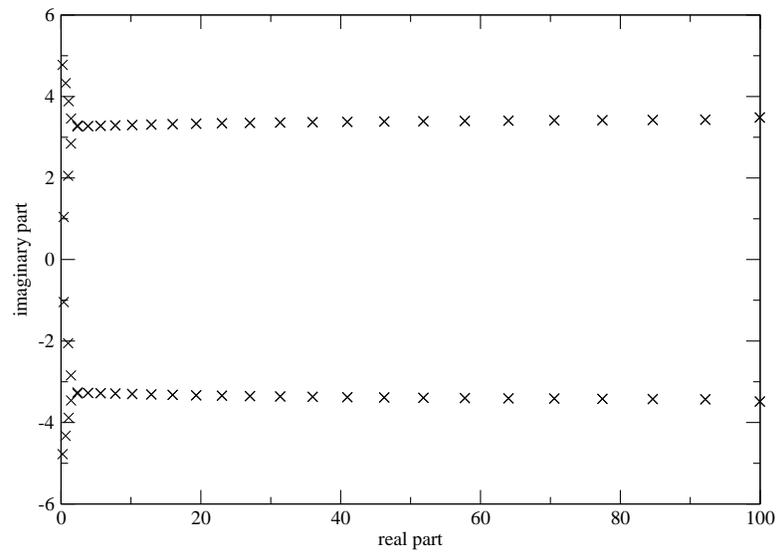}
\caption{ \label{Figure10} Spectrum of the matrix ${\mathcal  A}_{a,N}$ \Ref{eq110} (finite difference scheme) for $N=50$, $L=10$ and $a=1$.}
\end{figure}
\begin{figure}
  \centering
  \includegraphics[width=4in]{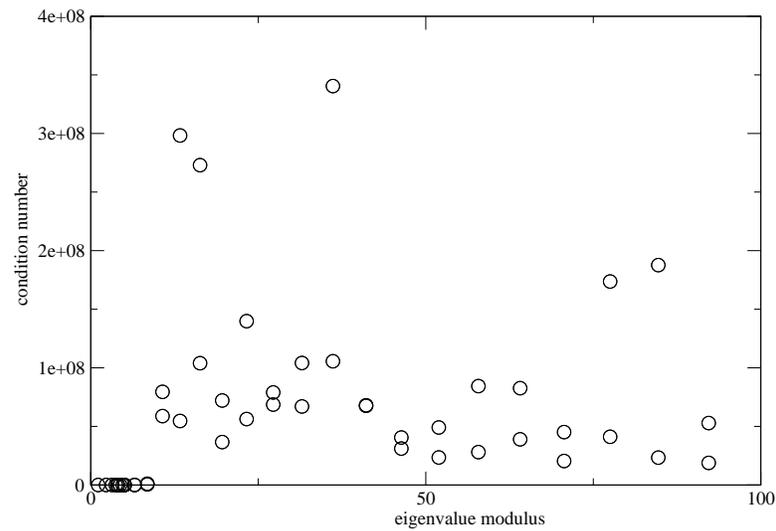}
\caption{ \label{Figure11} Condition number of the eigenvalues $\lambda_N$ in function of the modulus $\vert \lambda_N \vert $ for $N=50$, $L=10$ and $a=1$.} 
\end{figure}
\begin{figure}
  \centering 
  \includegraphics[width=4in]{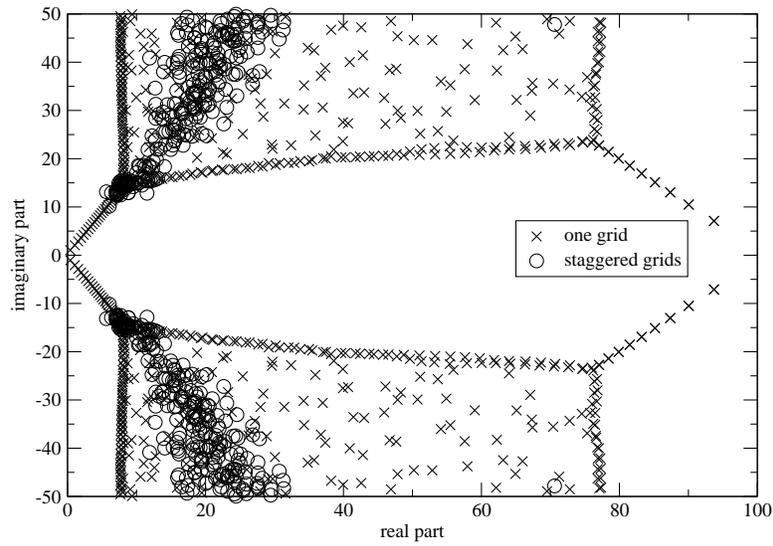}
\caption{ \label{Figure12} Average on the eigenvalues computed with the finite difference scheme using $11$ staggered grids for $a=1$, $N=1000$ and $L=10$. }
\end{figure}
\begin{figure}
  \centering
  \includegraphics[width=4in]{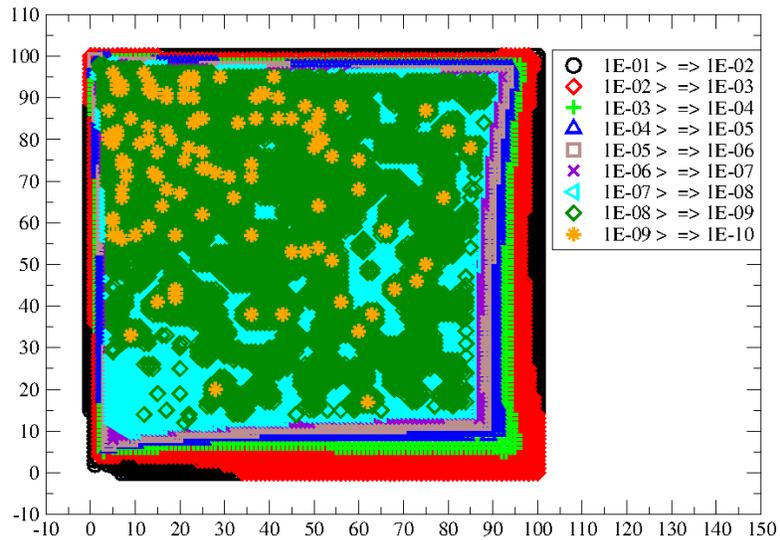}
\caption{ \label{Figure13} Computation of the pseudospectra \Ref{eq15} of the matrix ${\mathcal A}_{a,N}$ \Ref{eq110} (finite difference scheme), for $N=1000$, $L=10$ and $a=1$. }
\end{figure}
\begin{figure}
  \centering
  \includegraphics[width=4in]{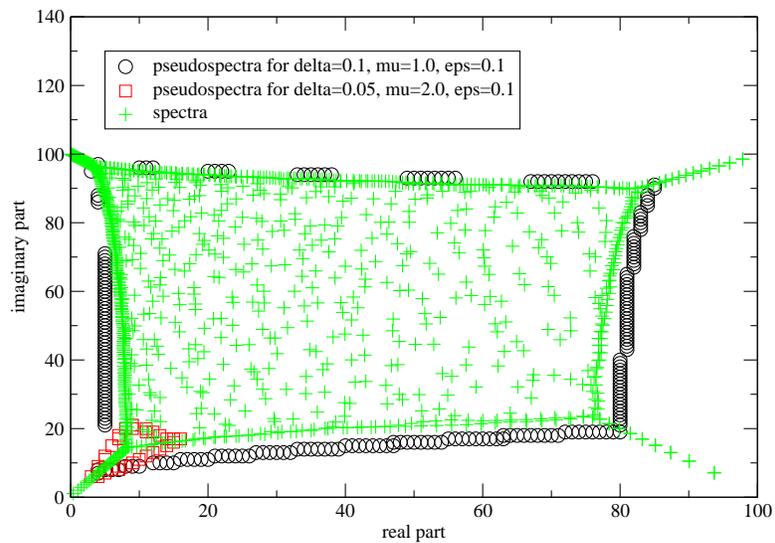}
\caption{ \label{Figure101} Computation of the pseudospectra \Ref{Pseudo3} (finite difference scheme), for $N=1000$, $L=10$ and $a=1$. }
\end{figure}
\newpage

\begin{figure}
  \centering
  \includegraphics[width=4in]{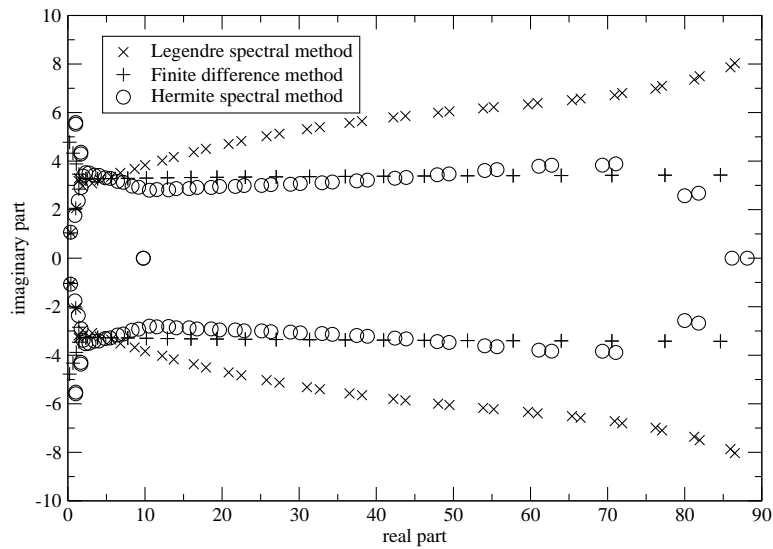} 
\caption{ \label{Figure14} Computation of the eigenvalues for $N=50$, $L=10$ and $a=1$ with the Legendre spectral method. Comparison with the spectral Hermite method and the finite difference method is done.}
\end{figure}


\begin{figure}
  \centering
  \includegraphics[width=5in]{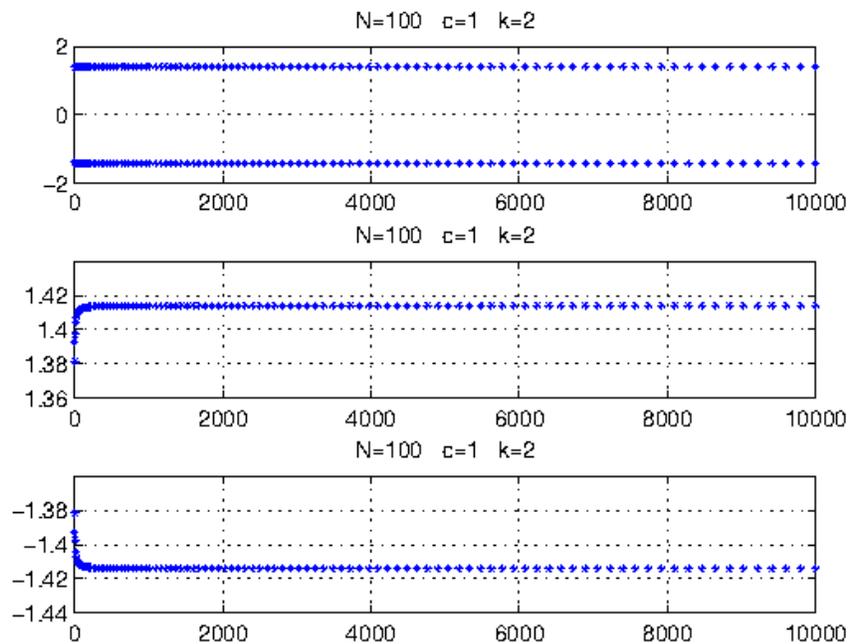} 
  \caption
   {\label{zomec1} Eigenvalues of the matrix $\mathcal{A}_c$ for $N=100$, $c=1$ and $k=2$. This figure represents a zoom for the case $c=1$.}\end{figure}
\newpage
\phantom{aaaa}
\begin{figure}
  \centering
  \includegraphics[width=6in]{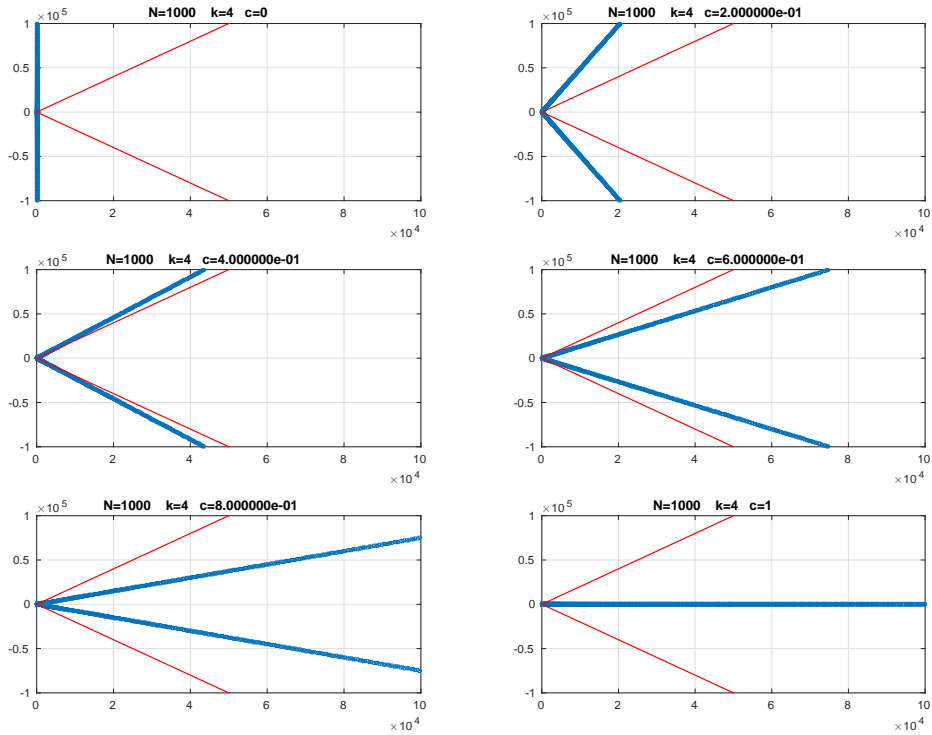} 
 \caption{\label{les6cassimplex4}  Eigenvalues of the matrix $\mathcal{A}_c$ for $N=1000$, $k=4$ and $c=0$, $0.2$, $0.4$, $0.6$, $0.8$, $1$. In the first three figures we can see the cases $c=0$, $0.2$, $0.4$. In the last three figures we can see the cases $c=0.6$, $0.8$, $1$.} \end{figure}
\begin{figure}
  \centering
  \includegraphics[width=4in]{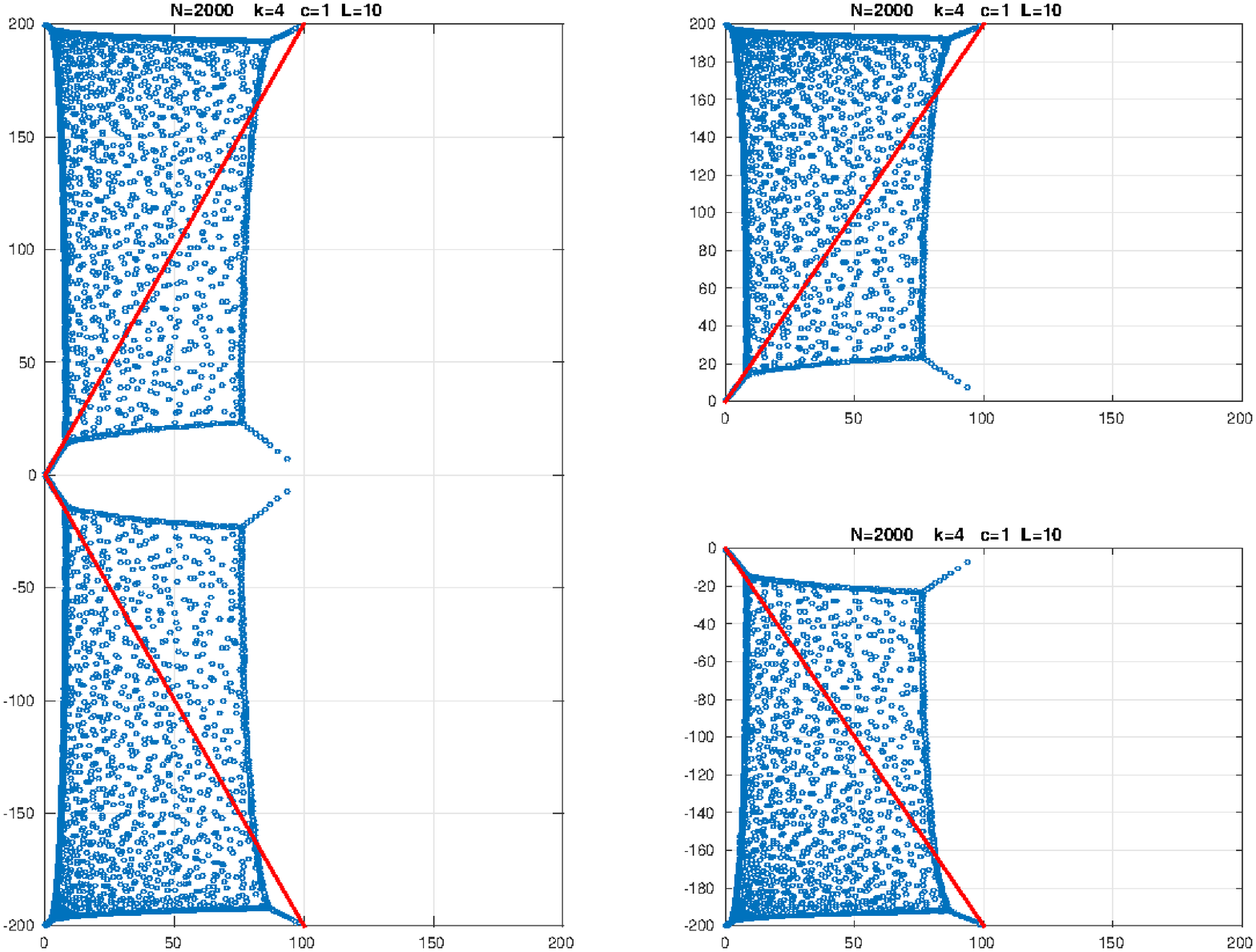} 
  \caption
   {\label{N2000k4c1L10tt3fig} Eigenvalues of the matrix $\mathcal{A}_c$ for $N=2000$, $k=4$, $L=10$, $c=1 $.}\end{figure}
\newpage
\begin{figure}
  \centering
  \includegraphics[width=4in]{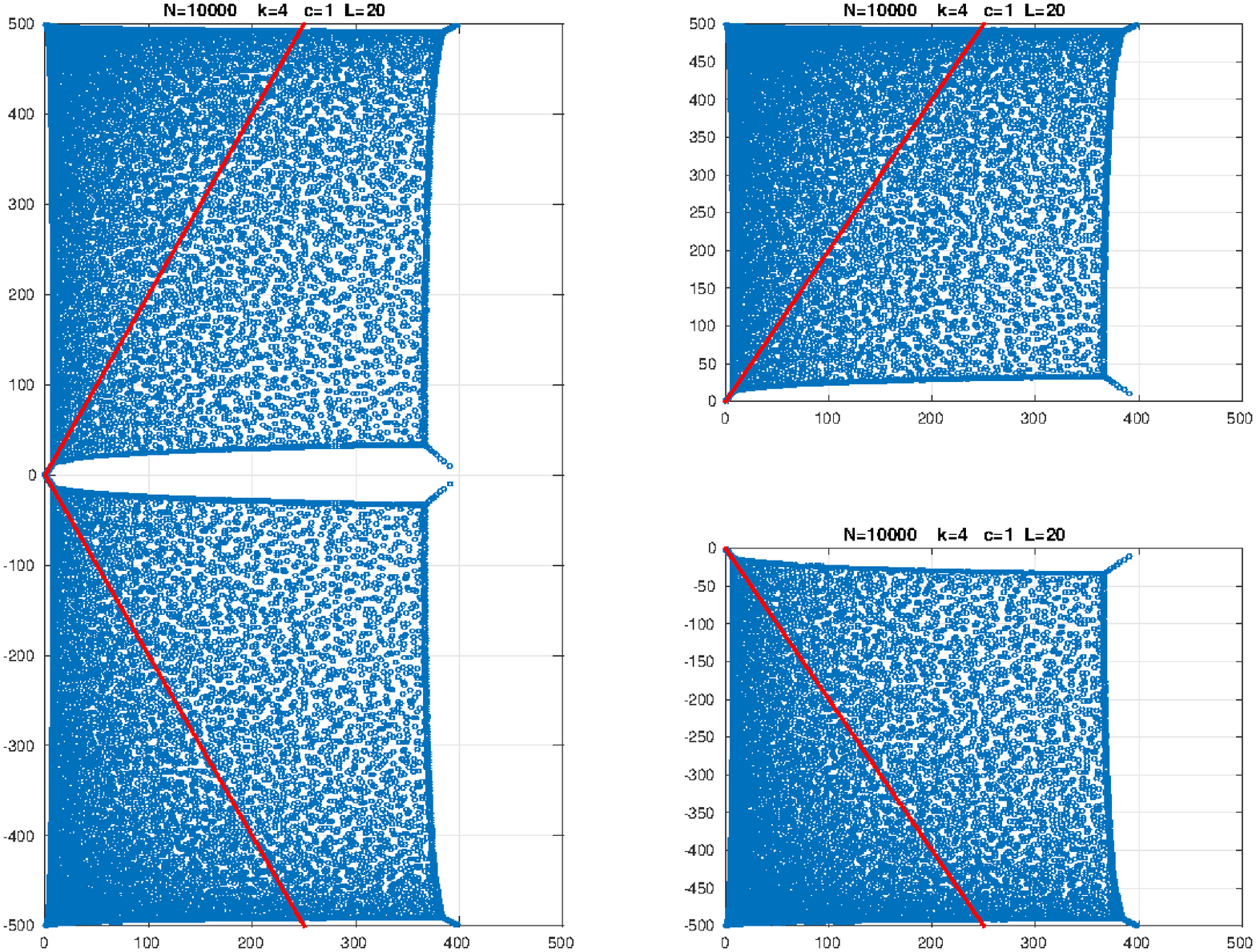} 
  \caption
   {\label{N10000k4c1L20tt3fig} Eigenvalues of the matrix $\mathcal{A}_c$ for $N=10000$, $k=4$, $L=20$, $c=1 $.}\end{figure}
\begin{figure}
  \centering
  \includegraphics[width=4in]{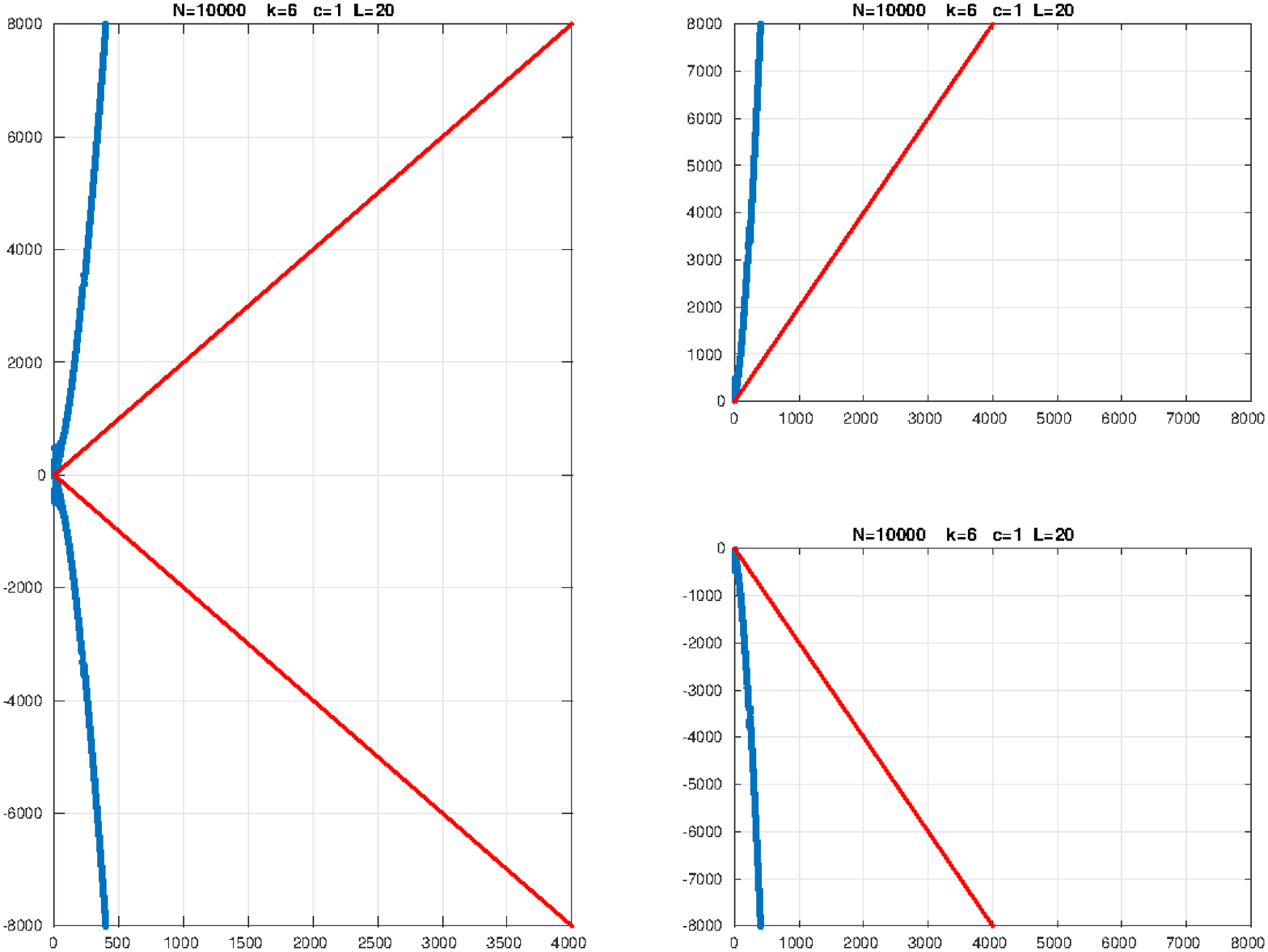} 
  \caption
   {\label{N10000k6c1L20tt3figs} Eigenvalues of the matrix $\mathcal{A}_c$ for $N=10000$, $k=6$, $L=20$, $c=1$. The figure on the left corresponds to $0\leq \Re \lambda \leq 10000$ and $-10000 \leq \Im \lambda \leq 10000$. The two figures on the right correspond, up to $0\leq \Re \lambda \leq 10000$ and $0 \leq \Im \lambda \leq 10000$, down to $0\leq \Re \lambda \leq 10000$ and $-10000 \leq \Im \lambda \leq 0$.}\end{figure}
\begin{figure}
  \centering
  \includegraphics[width=4in]{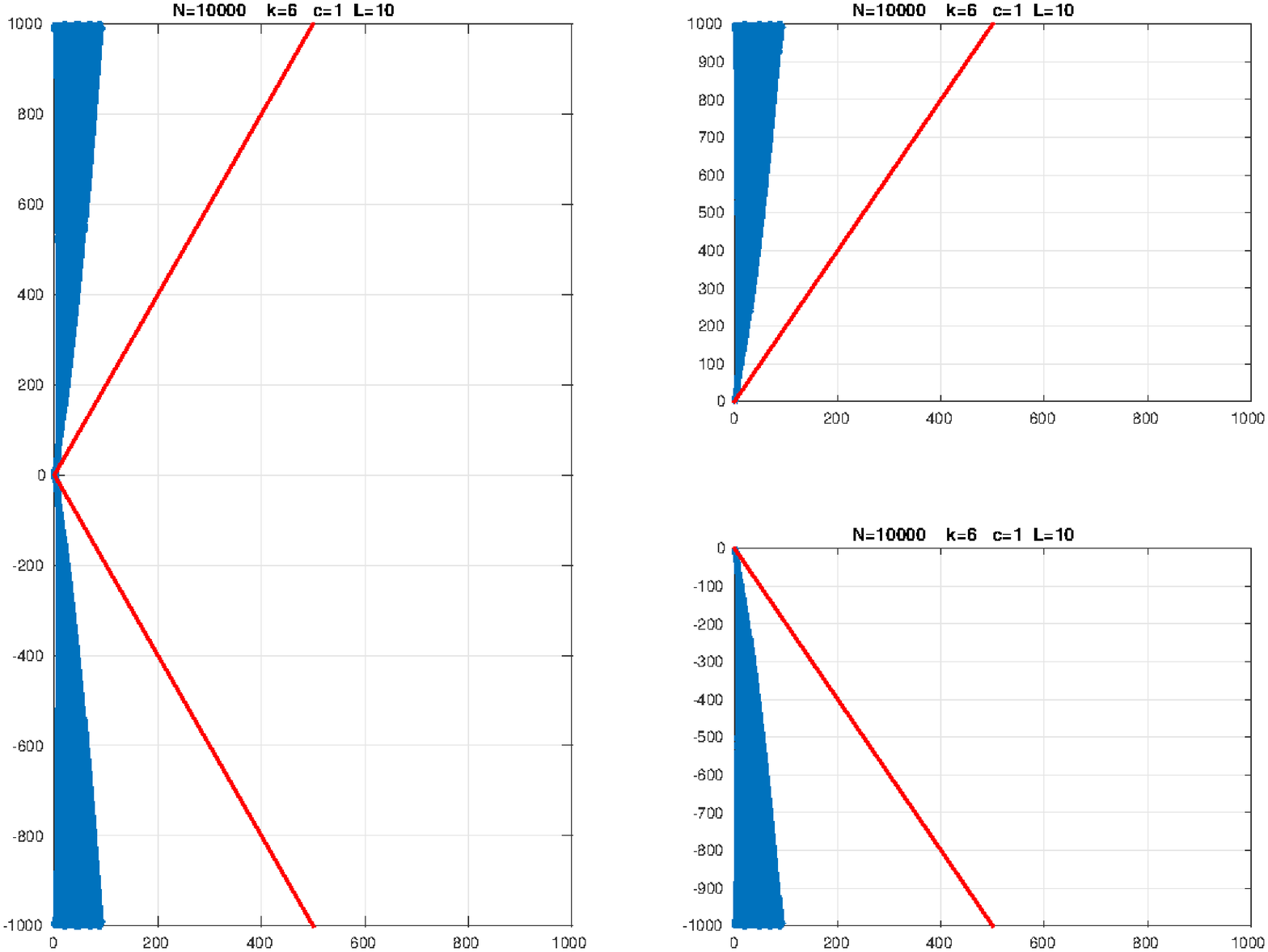} 
  \caption
   {\label{N10000k6c1L10tt3fig} Eigenvalues of the matrix $\mathcal{A}_c$ for $N=10000$, $k=6$, $L=10$, $c=1$. The figure on the left corresponds to $0\leq \Re \lambda \leq 1500$ and $-1500 \leq \Im \lambda \leq 1500$. The two figures on the right correspond, up to $0\leq \Re \lambda \leq 1500$ and $0 \leq \Im \lambda \leq 1500$, down to $0\leq \Re \lambda \leq 1500$ and $-1500 \leq \Im \lambda \leq 0$.}\end{figure}

\newpage

\phantom{aaaa}

\end{appendices}

\end{document}